\documentclass[12pt]{article}
\usepackage{amsfonts,latexsym}
\newtheorem{theorem}{Theorem}[section]
\newtheorem{lemma}[theorem]{Lemma}
\newtheorem{corollary}[theorem]{Corollary}
\newtheorem{proposition}[theorem]{Proposition}

\newcommand{\C}{{\mathbb C}}
\newcommand{\Z}{{\mathbb Z}}
\newcommand{\T}{{\mathbb T}}

\newcommand{\qed}{\hfill $\Box$}
\newcommand{\proof}{{\em Proof. }}
\newcommand{\be}{\begin{equation}}
\newcommand{\ee}{\end{equation}}
\newcommand{\bq}{\begin{eqnarray}}
\newcommand{\eq}{\end{eqnarray}}
\newcommand{\nn}{\nonumber}
\newcommand{\ba}{\begin{array}}
\newcommand{\ea}{\end{array}}
\newcommand{\ds}{\displaystyle}
\newcommand{\ts}{\textstyle}
\renewcommand{\Re}{\mbox{\rm Re\,}}

\renewcommand{\ll}{\ell_2}
\newcommand{\wt}[1]{\widetilde{#1}}
\newcommand{\iv}{^{-1}}
\newcommand{\iy}{\infty}
\newcommand{\ro}{\varrho}
\newcommand{\la}{\lambda}
\newcommand{\bt}{\beta}
\newcommand{\Ga}{\Gamma}
\newcommand{\diag}{\mbox{\rm diag\,}}
\newcommand{\ind}{\mbox{\rm ind\,}}
\newcommand{\tr}{\mbox{\rm tr\,}}
\newcommand{\slim}{\mbox{\rm s}-\lim}
\newcommand{\B}{B^1_1}
\newcommand{\GB}{G_0B^1_1}
\newcommand{\Ci}{C^\iy}
\newcommand{\PC}[1]{PC[#1]}
\newcommand{\PCa}[1]{PC_{\rm I}[#1]}
\newcommand{\PCb}[1]{PC_{\rm II}[#1]}
\newcommand{\PCc}[1]{PC_{\rm III}[#1]}

\renewcommand{\th}{\theta}

\newcommand{\lb}[2]{\makebox[#1][l]{$ #2 $}}

\newcommand{\ta}{\wt{a}}
\newcommand{\tb}{\wt{b}}
\newcommand{\ph}[1]{\phi^{(#1,\bt)}}
\newcommand{\hsp}{\hspace{\arraycolsep}} % spacing at `=' in display 
\newcommand{\eql}{\hsp=\hsp}

\begin{document}

\title{Asymptotic formulas for determinants of a sum of finite 
Toeplitz and Hankel matrices}
\author{Estelle L. Basor\thanks{ebasor@calpoly.edu. 
          Supported in part by NSF Grant DMS-9623278.}\\
               Department of Mathematics\\
               California Polytechnic State University\\
               San Luis Obispo, CA 93407, USA
        \and
        Torsten Ehrhardt\thanks{tehrhard@mathematik.tu-chemnitz.de.
          Supported in part by a DAAD Grant 213/402/537/5.}\\
       		Fakult\"{a}t f\"{u}r Mathematik\\
         	Technische Universit\"{a}t Chemnitz\\
\date{}        	09107 Chemnitz, Germany}
\maketitle              
 
\begin{abstract}
The purpose of this paper is to describe asymptotic formulas for 
determinants of a sum of finite Toeplitz and Hankel matrices with
singular generating functions. The formulas are similar to those
of the analogous problem for finite Toeplitz matrices for a certain 
class of symbols. However, the appearance of the Hankel matrices 
changes the nature of the asymptotics in some instances depending
on the location of the singularities. Several concrete examples
are also described in the paper.
\end{abstract} 

%%%%%%%%%%%%%%%%%%%%%%%%%%%%%%%%%%%%%%%%%%%%%%%%%%%%%%%%%%%%%%%%%%%

\section{Introduction}

In the theory of random matrices \cite{M} one is led naturally to
consider the asymptotics of determinants of Fredholm operators of
the form $I + W + H$ where $W$ is a finite Wiener--Hopf operator and
$H$ is a finite Hankel operator.  This problem arises when
investigating the probability distribution function of a random
variable thought of as a function of the eigenvalues of a positive
random Hermitian matrix.  The random matrix connections show that the
constant term in the asymptotic expansion of determinants is
fundamentally connected to the mean and variance of the distribution
function.  We will not describe the random matrix connections any
further, but simply refer the reader to \cite{M} for general
information and also \cite{F,B,BT,Be} for more specific tie--ins to
the random variable problem.

The focus of this paper is to study the discrete analogue of this
problem.  This is not exactly the desired situation for those
interested in random matrix theory.  However, it is a natural starting
place for cases where the random variable is discontinuous, since then
the discrete nature of the computations make things a bit more
accessible and the mathematical questions that arise are quite
interesting in themselves.

The discrete analogue of the Fredholm determinant problem is precisely
to find an asymptotic expansion of the determinants of the following
Toeplitz + Hankel matrices \bq M_n(\phi)&=&T_n(\phi)+H_n(\phi).  \eq
Here the $n\times n$ Toeplitz and Hankel matrices are defined as usual
by \be T_n(\phi) \eql \big(\phi_{j-k}\big)_{j,k=0}^{n-1},\qquad\quad
H_n(\phi) \eql \big(\phi_{j+k+1}\big)_{j,k=0}^{n-1}.  \ee The entries
$\phi_k$ are given by the $k$-th Fourier coefficient of $\phi$ where
$\phi$ is a sufficiently well-behaved function on the circle.  If
$\phi$ is continuous, even, and sufficiently smooth, then the
continuous analogue of the problem (i.e.~the Toeplitz + Hankel
matrices are replaced by finite Wiener--Hopf + Hankel operators) is
solved in \cite{B}.  There it is shown that the asymptotics are very
similar to the ones given in the Szeg\"{o}--Kac--Widom Strong Limit
Theorem.  Indeed, it is only in the constant, or third order term that
the answers differ.  This is no surprise since if $\phi$ is
continuous, then the Toeplitz operator is perturbed by a compact
Hankel operator.  However, if the symbol $\phi$ is singular, then the
problem, as in the Toeplitz case, is not easy to solve.

The purpose of this paper is to compute the asymptotics of $\det
M_n(\phi)$ as $n\to\iy$ for certain piecewise continuous functions
$\phi$.  The main general result that we will obtain is as follows.
We consider piecewise continuous functions of the form
\bq\label{F.M}
\phi(e^{i\th}) &=& b(e^{i\th}) t_{\bt_+}(e^{i\th})
t_{\bt_-}(e^{i(\th-\pi)}) \prod_{r=1}^R 
t_{\bt_r}(e^{i(\th-\th_r)}),\label{1}
\eq
where
\bq 
t_{\bt}(e^{i\th}) &=& \exp(i\bt(\th-\pi)),\qquad\quad 
0<\th<2\pi,\label{1a}
\eq
and $b$ is a smooth nonvanishing function defined on the circle with
winding number zero.  We also need conditions on the parameters
$\bt_+,\bt_-,\bt_1,\dots,\bt_R$.  These conditions on the ``size'' of
the jumps and the precise smoothness conditions on $b$ will be
described later on.  In addition, we have to assume that
$\th_1,\dots,\th_R\in(-\pi,0)\cup(0,\pi)$ are certain distinct
numbers satisfying $\th_r+\th_s\neq0$ for each $r$ and $s$.
The latter condition excludes piecewise continuous functions
with jumps at both a point on the unit circle and its complex
conjugate. However, the function $\phi$ may have jumps at the points
$1$ and $-1$. Our conditions on $b$ guarantee that the following
functions   
\bq
b_+(t) &=& \exp\left(\,\sum_{k=1}^\iy t^k[\log b]_k\right),
\qquad t\in\T,\label{2}\\
b_-(t) &=& \exp\left(\,\sum_{k=1}^\iy t^{-k}[\log b]_{-k}\right),
\qquad t\in\T,\label{3}
\eq
are well defined and smooth.  Here $[\log b]_k$ denotes the $k$-th
Fourier coefficient of the logarithm of $b$.
Then the asymptotic formula reads
\bq
 \det M_n(\phi) &\sim& G[b]^{n}n^{\Omega_M}E_M\label{4}
\eq
as $n\to\iy$, where
\bq
G[b] &=& \exp[\log b]_{0},\label{5}\\ 
\Omega_M &=& -\frac{3\bt_+^2}{2}-\frac{\bt_+}{2}
  -\frac{3\bt_-^2}{2}+\frac{\bt_-}{2}
  -\sum_{r=1}^{R}\bt_{r}^{2},\label{6}\\
E_M &=& E[b] F[b] \nn\\
 &&\times\; G(1+\bt_+)G(1-\bt_+)G(1/2-\bt_+)G(1/2)\iv
  (2\pi)^{\bt_+/2}\, 2^{3\bt_+^2/2}\nn\\ 
 &&\times\; G(1+\bt_-)G(1-\bt_-)G(3/2-\bt_-)G(3/2)\iv
  (2\pi)^{\bt_-/2}\, 2^{3\bt_-^2/2}\nn\\ 
 &&\times\prod_{r=1}^R G(1+\bt_{r})G(1-\bt_{r})
  \left(1-e^{-i\th_r}\right)^{\bt_r^2/2+\bt_r/2} 
  \left(1+e^{-i\th_r}\right)^{\bt_r^2/2-\bt_r/2}\nn\\
 &&\times\; b_+(1)^{2\bt_+} b_-(1)^{-\bt_+}
  b_+(-1)^{2\bt_-} b_-(-1)^{-\bt_-}\, 2^{3\bt_+\bt_-}\nn\\
 &&\times\prod_{r=1}^R b_+(e^{i\th_r})^{\bt_r}
  b_-(e^{i\th_r})^{-\bt_r} b_+(e^{-i\th_r})^{\bt_r}\nn\\
 &&\times\prod_{r=1}^R
   \left(1-e^{-i\th_r}\right)^{2\bt_+\bt_r}
   \left(1-e^{i\th_r}\right)^{\bt_+\bt_r}
   \left(1+e^{-i\th_r}\right)^{2\bt_-\bt_r}
   \left(1+e^{i\th_r}\right)^{\bt_-\bt_r}\nn\\
 &&\times\prod_{1\leq s<r\leq R}
   \left(1-e^{i(\th_s-\th_r)}\right)^{\bt_r\bt_s}
   \left(1-e^{i(\th_r-\th_s)}\right)^{\bt_r\bt_s}
   \left(1-e^{-i(\th_s+\th_r)}\right)^{\bt_r\bt_s}.\label{7}
\eq
The constants $E[b]$ and $F[b]$ are defined by
\bq
E[b] &=&\exp\left(\,\sum_{k=1}^\iy k[\log b]_k[\log b]_{-k}\right),
 \label{8}\\
F[b] &=&\left(\frac{b_+(1)}{b_+(-1)}\right)^{1/2}\,
 \exp\left(-\frac{1}{2}\sum_{k=1}^\iy k[\log b]_k^2\right),\label{8a}
\eq		
and $G(\ast)$ is the Barnes G--function \cite{Bar,WW} defined by
\bq
 G(1+z) &=&(2\pi)^{z/2}e^{-(z+1)z/2-\gamma_{E}z^{2}/2}
 \prod_{k=1}^\iy\left((1+z/k)^{k}e^{-z+z^{2}/2k}\right)\label{9}
\eq
with $\gamma_{E}$ being Euler's constant.

It is interesting to compare this asymptotic formula with the
corresponding formula for Toeplitz determinants.  The asymptotic
expansion of Toeplitz determinants for a certain class of singular
generating functions is described by the Fisher--Hartwig conjecture.
For instance, it is well known \cite{B79,Bo} that if $\phi$ is of the
form
\bq\label{F.T}
\phi(e^{i\th}) &=& b(e^{i\th}) \prod_{r=1}^R 
t_{\bt_r}(e^{i(\th-\th_r)}),
\eq
where $b$ is a sufficiently smooth and nonvanishing on the unit
circle, $\th_1,\dots,\th_R\in(-\pi,\pi]$ are distinct numbers, and
if $|\Re \bt_r|<1/2$ holds for each $1\le r\le R$, then the
asymptotic behavior of the Toeplitz determinants is given by
\bq\label{A.T}
 \det T_n(\phi) &\sim& G[b]^nn^{\Omega_T} E_T
\eq
as $n\to\iy$, where 
\bq
\Omega_T &=&
  -\sum_{r=1}^{R}\bt_{r}^{2},\\
E_T &=& E[b]\, \prod_{r=1}^R G(1+\bt_r)G(1-\bt_r)\nn\\
  &&\times \prod_{r=1}^R
    b_+(e^{i\th_r})^{\bt_r} b_-(e^{i\th_r})^{-\bt_r}
  \;\prod_{1\le r\neq s\le R}
  \Big(1-e^{i(\th_s-\th_r)}\Big)^{\bt_r\bt_s}.
\eq
A general account of the Fisher--Hartwig conjecture can be found in
\cite{BS} and more recent work in \cite{BT1,ES,Eh}.  It has been proved in
many cases and reformulated in others.

The paper is organized as follows.  In Section \ref{s.smooth} we
compute the asymptotics of $\det M_n(\phi)$ in the case of smooth
nonvanishing functions $\phi$ with winding number zero.  As in the
continuous analogue, the asymptotics are very similar to the
Toeplitz case.  In fact, we prove that the quotient $\det
M_n(\phi)/\det T_n(\phi)$ converges to a nonzero constant.

In Section \ref{s.operator} we recall several operator
theoretic results, in particular those related to
Toeplitz operators and Toeplitz + Hankel operators.

In Section \ref{s.limit} we prove the asymptotic formula
(\ref{4}) in the special case of piecewise continuous functions
(\ref{F.M})
without jumps at $1$ and $-1$ (i.e.~with $\bt_+=\bt_-=0$), without
jumps at both a point on the unit circle and its complex conjugate
(i.e.~$\th_r+\th_s\neq0$) and under the assumption
$|\Re \beta_r|<1/2$ , $1\le r\le R$.
As before we show that $\det M_n(\phi)/\det T_n(\phi)$
converges to a nonzero constant (althought the corresponding Hankel
operator is not compact).  
In other words, if the symbol has jumps in the ``proper'' locations,
then the asymptotic behavior is again like in the Toeplitz case.
Note that the condition on the location of the jumps is extremely
important in the Toeplitz + Hankel case, as contrasted to the Toeplitz
case.

In Section \ref{s.local} we show that the quotient
$\det M_n(\phi\psi)/(\det M_n(\phi)\det M_n(\psi)$
converges to a nonzero constant for certain functions
$\phi$ and $\psi$.  This result allows
us to localize at certain points on the unit circle, in particular at
$1$ or $-1$.  However, it is not possible to localize at a point on
the unit circle and its complex conjugate. Thus the localization result
reduces the asymptotics for general piecewise continuous functions to
those for functions $t_\bt(e^{i\th})$ and $t_\bt(e^{i(\th-\pi)})$
with a single pure jump at $1$ or $-1$ and to those for functions
$t_{\bt_r^+}(e^{i(\th-\th_r)})t_{\bt_r^-}(e^{i(\th+\th_r)})$
with two pure jumps. Note that it is exactly this last class of
functions for which we are not able to determine the asymptotics
in general.

In Section \ref{s.1jump} we then consider the case of functions
$t_\bt(e^{i\th})$ and $t_\bt(e^{i(\th-\pi)})$.  In this case, the jump
at $1$ or $-1$ on the circle changes the nature of the asymptotics in
the second order term.  The computations are based on the fact that
the corresponding matrices are Cauchy matrices times certain diagonal
matrices.  Finally, this result in conjunction with the localization
and the results of Section \ref{s.limit} gives the above mentioned
main general result (\ref{4}).

In the last section, we illustrate with additional concrete
examples.  The first class of examples is for piecewise continuous
functions with two jumps either at $\pm1$ or at $\pm i$.  These
functions are special cases where the parameters describing the jumps
are connected with each other in some way.  The significance is that
they show that the one jump results, for functions with jumps at $i$
and $-i$ cannot be pieced together to obtain the asymptotics for a
symbol that has jumps at both these points.  It should be pointed out
that this does work for the points $1$ and $-1$ and this is confirmed
by the examples.  Note that one special case of these
examples is a piecewise constant even function with jumps at
$\pm i$.

The second class of examples in the last section is for the even
functions
\be
u_\alpha(e^{i\th}) \eql (2-2 \cos\th)^{\alpha}
\quad\mbox{ and }\quad
u_\alpha(e^{i(\th-\pi)}) \eql (2+2 \cos\th)^{\alpha}.
\ee
These function are singular or zero at $1$ or $-1$, respectively, and
are particularly interesting since they represent
a more general class of examples of even functions.

For random matrix theory even functions are most important.  It would
be helpful in the future to extend these results to that case and also
to the continuous analogue of Wiener--Hopf + Hankel operators.
But we believe the present paper is a good start and
leave the other questions to some other time.

\section{Operator theoretic results and Toeplitz + Hankel
determinants in the case of smooth functions}
\label{s.smooth}

In this section, we compute the asymptotic behavior of determinants
of Toeplitz + Hankel matrices $M_n(\phi)$ in the case of smooth
nonvanishing functions defined on the unit circle with winding number
zero.  What we exactly mean by smoothness will be explained shortly.
In the first part of this section, however, we will recall
certain operator theoretic results that will be needed later on.

We first introduce the following linear bounded operators acting
on the Hilbert space $\ll=\ll(\Z_+)$ of one-sided square-summable
sequences.  Given $\phi\in L^\iy(\T)$, define
\bq
M(\phi) &=& T(\phi)+H(\phi)
\eq
where the Toeplitz and Hankel operators are given by
the infinite matrices
\be
T(\phi) \eql \big(\phi_{j-k}\big)_{j,k=0}^\iy,\qquad\quad
H(\phi) \eql \big(\phi_{j+k+1}\big)_{j,k=0}^\iy.
\ee
Note that the Hardy spaces $H^\iy$ and $\overline{H^\iy}$ consist of
those functions $\phi\in L^\iy(\T)$ for which the Fourier coefficients
$\phi_n$ vanish for each $n<0$ or $n>0$, respectively.  We also write
\bq
\wt{\phi}(e^{i\th}) &=& \phi(e^{-i\th}),
\eq
and call $\phi$ even if $\wt{\phi}=\phi$.
Finally, we introduce the projection $P_n$ acting on $\ll$ by
\bq
P_n (f_{0},f_{1},\dots)  &=& 
(f_0,f_1,\dots,f_{n-2},f_{n-1},0,0,\dots).
\eq
Note that $T_n(\phi)=P_nT(\phi)P_n$, $H_n(\phi)=P_nH(\phi)P_n$
and $M_n(\phi)=P_nM(\phi)P_n$.

It is well known that Toeplitz and Hankel operators are related
to each other by the formulas
\bq
T(\phi\psi) &=& T(\phi)T(\psi)+H(\phi)H(\wt{\psi}),\label{f0.Txx}\\
H(\phi\psi) &=& T(\phi)H(\psi)+H(\phi)T(\wt{\psi}).\label{f0.Hxx}
\eq
If $\psi_+\in H^\iy$ and $\psi_-\in \overline{H^\iy}$, then we have
\bq
T(\psi_-\phi\psi_+)      &=& T(\psi_-)T(\phi)T(\psi_+),\label{f0.Txxx}\\
H(\psi_-\phi\wt{\psi}_+) &=& T(\psi_-)H(\phi)T(\psi_+).\label{f0.Hxxx}
\eq
Combining equations (\ref{f0.Txx}) and (\ref{f0.Hxx}),
it follows that
\bq\label{f0.TM}
M(\phi\psi) &=& T(\phi)M(\psi)+H(\phi)M(\wt{\psi}).
\eq
This implies 
\bq\label{f0.Mxx}
M(\phi\psi) &=& M(\phi)M(\psi) + H(\phi)M(\wt{\psi}-\psi).
\eq
If $\psi$ is even, then the latter equation simplifies to
\bq\label{f0.Mx}
M(\phi\psi) &=& M(\phi)M(\psi).
\eq
These and other results concerning $M(\phi)$ are discussed and
proved in \cite{BaEh}.

Two important notions are stability and strong convergence.
Let $A_n$ be a sequence of operators.  This sequence is said
to be stable if there exists an $n_0$ such that the operators $A_n$
are invertible for each $n\geq n_0$ and ${\sup}_{n\geq n_0}
\|A_{n}^{-1}\|<\infty$.  Moreover, we say that $A_n$ converges
strongly on $\ll$ to an operator $A$ as $n\to\iy$ if $A_nx\to Ax$ in
the norm of $\ll$ for each $x\in\ll$.  When dealing with finite
matrices $A_n$, we identify the matrices and their inverses with
operators acting on $\ll$.  It is interesting to note that
stability is related to strong convergence of the inverses
(and their adjoints) in the following sense.

\begin{lemma}\label{l1.1}
Suppose that $A_n$ is a stable sequence such that $A_n\to A$ and
$A_n^*\to A^*$ strongly.  Then $A$ is invertible, and
$A_n\iv\to A\iv$ and $(A_n\iv)^*\to (A\iv)^*$ strongly.
\end{lemma}

Another important set of operators is the ideal of trace class operators (see e.g.
\cite{GK}).  For such operators, the trace ``$\tr A$'' and the
operator determinant ``$\det(I+A)$'' are well defined and continuous
with respect to $A$ in the trace class norm.  The following result
shows the connection with strong convergence.

\begin{lemma}\label{l1.1b}
Let $B$ be a trace class operator and suppose that $A_n$ and $C_n$
are sequences such that $A_n\to A$ and $C_n^*\to C^*$ strongly. Then 
$A_nBC_n\to ABC$ in the trace class norm.
\end{lemma}

We proceed with describing the smoothness conditions. We therefore introduce
certain function spaces. Let $F\ll^\bt$ stand for
the Banach space of all functions $b\in L^1(\T)$ for which
\bq
\|b\|_{F\ll^\bt} &:=& \left(\sum_{n=-\iy}^\iy(1+|n|)^{2\beta}|b_n|^2
\right)^{1/2}\;<\;\iy,
\eq
and let $W$ denote the Wiener algebra.  It is well known that
$F\ll^{1/2}\cap W$ is a Banach algebra of continuous functions on the
unit circle.  The Besov class $\B$ is the Banach algebra of all
functions $b\in L^1(\T)$ for which
\bq
\|b\|_{\B} &:=&
\int_{-\pi}^\pi\frac{1}{y^2}\int_{-\pi}^\pi\left|\,
b(e^{ix+iy})+b(e^{ix-iy})-2b(e^{ix})\,\right|\,dx\,dy
\;<\;\iy.
\eq
Using results of Peller \cite{Pe} one can show that $b\in\B$ if and
only if both $H(b)$ and $H(\wt{b})$ are trace class operators.
Moreover, the Riesz projection acts boundedly on $\B$. An
equivalent norm in $\B$ is given by $|b_0|+\|H(b)\|_{1}+\|H(\wt{b})\|_{1}$.
Finally, one has the following continuous and dense embeddings
\bq
F\ll^{\beta}  \hsp\subset\hsp  \B  \hsp\subset\hsp
\Big(F\ll^{1/2}\cap W\Big)  \qquad\quad\mbox{if }\beta>1.
\eq
For more information on these and related classes of smooth functions
we refer the reader to \cite{BS} and the literature cited there.

The Besov class $\B$ exactly fits our purposes in the sense that
the function $b$ appearing in (\ref{F.M}) will be assumed to be in $\B$.
In order to simplify notation we denote by $\GB$ the group of all
nonvanishing functions in $\B$ with winding number zero.  Remark that
the asymptotic formula for Toeplitz determinants as given in
(\ref{A.T}) has been proved for functions of the form (\ref{F.T}) with
$b\in\GB$ (see e.g.  \cite{BS}).  The following proposition shows that 
all definitions involving the function $b$ which were made in the 
introduction make sense.
\begin{proposition}
Let $b\in\GB$. Then $b$ possesses a logarithm $\log b\in\B$, and
the constants $G[b]$, $E[b]$ and $F[b]$ as well as the functions
$b_+$ and $b_-$ make sense. Moreover,
\bq
b(e^{i\th}) &=& b_-(e^{i\th})G[b]b_+(e^{i\th}),\qquad 0\le\th<2\pi,\nn
\eq
is the (normalized) Wiener--Hopf factorization with
$b_+\in G(\B\cap H^\iy)$ and $b_-\in G(\B\cap\overline{H^\iy})$.
\end{proposition}
\proof
Approximating the function $b$ by polynomials and using that $\B$ is a
Banach algebra, one can show that $b$ possesses a logarithm $\log
b\in\B$.  Now one need only use the boundedness of the Riesz
projection and the fact that $\B$ is contained in $F\ll^{1/2}\cap W$.
\qed

We now establish a limit relation for the quotient of two determinants
with smooth generating functions.  In Section \ref{s.limit} this
relation will be generalized to the case of certain piecewise
continuous functions.
\begin{proposition}
\label{p0.2}
Let $b\in\GB$.  Then $\det M_n(b)/\det T_n(b)\to F(b)$ as $n\to\iy$
where $F(b)$ is the operator determinant defined by
\bq\label{f0.F}
F(b) &:=& \det\Big(I+T\iv(\phi)H(\phi)\Big).
\eq
\end{proposition}
\proof
It is well known \cite{GF} that under the above assumptions on $b$,
the sequence $T_n(b)$ is stable and the inverses converge strongly
on $\ll$ to the inverse of $T(b)$. Since $H(b)$ is trace class,
we obtain that
$$
\frac{\det M_n(b)}{\det T_n(b)} \eql
\det T_n\iv(b) M_n(b) \eql
\det\Big(P_n+T_n\iv(b)P_nH(b)P_n\Big)
$$
converges to the operator determinant $F(b)$.
Note that $P_n^*=P_n\to I$ strongly.
\qed

Next we establish the one main result of this section, the evaluation
of the operator determinant $F(b)$.  The computation is especially
remarkable since it relies on a ``differentiation'' trick, which was
recently used in a modified form in \cite{B}.

\begin{theorem}\label{t0.2}
Let $b\in\GB$. Then $F(b)=F[b]$.
\end{theorem}
\proof
Using the Wiener--Hopf factorization of $b$ and formula
(\ref{f0.Txxx}) it is easily seen that $T(b)=G[b]T(b_-)T(b_+)$, and
hence the inverse equals
\bq
T\iv(b) &=& G[b]\iv T(b_+\iv)T(b_-\iv).\nn
\eq
{From} equation (\ref{f0.Hxxx}) it follows that
\bq
T\iv(b)H(b) &=& T(b_+\iv)T(b_-\iv)H(b_-b_+)\eql
T\iv(b_+)H(b_+).\label{f0.TiH}
\eq
Hence $F(b)=F(b_+)$.  Now let $c=\wt{b}_+b_+$.  One can conclude
analogously that $F(c)=F(b_+)$.  Thus we are left with the evaluation
of $F(c)$.  Obviously, $[\log c]_n=[\log b]_{|n|}$ for $n\in
\Z\setminus\{0\}$, and hence $\log\wt{c}=\log c$.  Since the Riesz
projection is bounded on $\B$, we obtain that $\log c\in\B$.  Next we
define the $\B$--valued analytic (in $\la$) function
\bq
c_\la &=& \exp(\la\log c),\qquad\qquad \la\in\C.\nn
\eq
Note that the derivative of $c_\la$ with respect to $\la$ equals 
$c_\la\log c$.
Let $Y(\la)$ be the analytic operator--valued function
\bq
Y(\la) &=& I+T\iv(c_\la)H(c_\la)\hsp=\hsp T\iv(c_\la)M(c_\la).\nn
\eq
It is easy to compute the inverse and the derivative of $Y(\la)$. We 
obtain
\bq
Y'(\la)Y\iv(\la) &=& -\,T\iv(c_\la)T(c_\la\log c)T\iv(c_\la)M(c_\la)
M\iv(c_\la)T(c_\la)\nn\\
&& +\,T\iv(c_\la)M(c_\la\log c)M\iv(c_\la)T(c_\la)\nn\\
&=& -\,T\iv(c_\la)T(c_\la\log c)+T\iv(c_\la)M(\log c)T(c_\la).\nn
\eq
Note that $\wt{c}_\la=c_\la$, and thus (\ref{f0.Mx}) implies
$M(c_\la\log c)=M(\log c)M(c_\la)$.  Because $\det Y(\la)\neq0$ for
all $\la$, the scalar function $y(\la)=\log\det Y(\la)$ is an entire
analytic function.  We conclude
\bq
\frac{dy}{d\la} &=& \frac{(\det Y)'}{\det Y}
\hsp=\hsp \tr Y'Y\iv\nn\\
&=& \tr\Big(-T(c_\la\log c)T\iv(c_\la)+M(\log c)\Big).\nn
\eq
Differentiating again yields
\bq
\frac{d^2y}{d\la^2} &=& \tr\Big(-T(c_\la\log^2c)T\iv(c_\la)+
T(c_\la\log c)T\iv(c_\la)T(c_\la\log c)T\iv(c_\la)\Big)\nn\\
&=& \tr\Big(-T(\log^2c)+T(\log c)T(\log c)\Big)
\hsp=\hsp-\,\tr H(\log c)H(\log \wt{c}).\nn
\eq
The last equality follows from
$T\iv(c_\la)=T(c_{\la,+}\iv)T(c_{\la,-}\iv)$ where
$c_\la=c_{\la,-}c_{\la,+}$ is the Wiener--Hopf factorization of
$c_\la$.  By repeated application of (\ref{f0.Txxx}), all occurring
functions $c_{\la,\pm}$ cancel each other.  Hence the second
derivative of $y$ does not depend on $\la$.  Note that
\bq
 y(0)  &=& \log\det I=0,\nn\\
 y'(0) &=& \tr\Big(-T(\log c)+M(\log c)\Big)
 \hsp=\hsp \tr H(\log c)\nn
\eq
since $c_\la|_{\la=0}\equiv 1$.  It follows that
\bq
y(\la) &=& -\,\la^2/2\;\tr H(\log c)H(\log \wt{c})+\la\;\tr H(\log 
c).\nn
\eq
Because $F(c)=\det Y(1) = \exp y(1)$, we obtain
\bq
F(c) &=&
\exp\Big(-1/2\,\tr H(\log c)H(\log \wt{c})+\tr H(\log c)\Big)\nn\\
&=& \exp\Big(-1/2\sum_{n=1}^\iy n [\log c]_n[\log c]_{-n}
+1/2\sum_{n=1}^\iy\Big\{[\log c]_n - (-1)^n [\log c]_n\Big\}\Big).\nn
\eq
Writing $\log c$ in terms of $\log b$, we arrive at
\bq
F(b) &=& \exp\Big(-1/2\,\sum_{n=1}^\iy n[\log b]_n^2
+1/2\,\log b_+(1)-1/2\,\log b_+(-1)\Big).\nn
\eq
This immediately implies the assertion.
\qed

Finally, we can combine the previous results with the well known
Szeg\"o--Widom Limit Theorem \cite{W}.  This limit theorem says that
$\det T_n(b)\sim G[b]^nE[b]$ as $n\to\iy$ for nonvanishing functions
$b\in F\ll^{1/2}\cap W$ with winding number zero.
\begin{corollary}
Let $b\in\GB$. Then
\bq
\det M_n(b) &\sim& G[b]^nE[b]F[b]
\qquad\mbox{ as }n\to\iy.\nn
\eq
\end{corollary}

\section{Further operator theoretic results}
\label{s.operator}

The proofs of the results that will presented in the following two
sections require further operator theoretic preliminaries.  In
particular, we need some results about Toeplitz operators and Hankel
operators as well as about Toeplitz + Hankel operators $M(\phi)$ (see
\cite{BS} and \cite{BaEh} for the general theory).
First of all, in addition to the projection $P_n$,
we define $Q_n=I-P_n$ and
\bq
W_n (f_{0},f_{1},\dots)  &=& (f_{n-1},f_{n-2},\dots 
,f_{1},f_{0},0,0,\dots),\nn\\  
V_n (f_{0},f_{1},\dots)  &=& 
(0,0,\dots,0,0,f_{0},f_{1},f_{2},\dots),\nn\\
V_{-n}(f_{0},f_{1},\dots)&=& (f_{n},f_{n+1},f_{n+2},\dots).\nn
\eq
It is easily seen that
$W_n^2=P_n$, $W_n=W_nP_n=P_nW_n$, $V_nV_{-n}=Q_n$ and $V_{-n}V_n=I$.
Note also that
\be
P_nT(\phi)V_n \eql W_nH(\wt{\phi}), \qquad\qquad
V_{-n}T(\phi)P_n \eql H(\phi)W_n. \label{f1.VW}\\
\ee
Moreover, we have
\bq
V_{-n}H(\phi) &=& H(\phi)V_n,\label{f1.VV}\\
W_nT_n(\phi)W_n &=& T_n(\wt{\phi}).\label{f1.WW}
\eq
Using equations (\ref{f1.VW}) we can write
\bq
P_nT(\phi)Q_nT(\psi)P_n &=& P_nT(\phi)V_nV_{-n}T(\psi)P_n \eql
W_nH(\wt{\phi})H(\psi)W_n.
\eq
Taking (\ref{f0.Txx}) into account we arrive at the
fundamental identity due to Widom \cite{W}
\bq\label{f1.Tn}
T_n(\phi\psi) &=& T_n(\phi)T_n(\psi)+P_nH(\phi)H(\wt{\psi})P_n+
W_nH(\wt{\phi})H(\psi)W_n.
\eq

The following result deals with strong convergence.
For brevity of notation, we will henceforth write
$A=\slim A_n$ if both $A_n\to A$ and $A_n^*\to A^*$ strongly.
\begin{proposition}
Let $\phi\in L^\iy(\T)$. Then
$$
\ba{rcl@{\qquad\qquad}rcl}
T(\phi) &=& \slim\; T_n(\phi), &
T(\wt{\phi}) &=& \slim\; W_nT_n(\phi)W_n,\\
M(\phi) &=& \slim\; M_n(\phi), &
T(\wt{\phi}) &=& \slim\; W_nM_n(\phi)W_n.
\ea
$$
\end{proposition}
\proof
The relations without the $W_n$'s are trivial because $P_n^*=P_n\to I$
strongly. Also the second assertion is easy to show by taking
account of (\ref{f1.WW}). Finally, using (\ref{f1.VW})
we obtain
$$
W_nM_n(\phi)W_n \eql T_n(\wt{\phi})+W_nV_{-n}T(\phi)P_n
\eql  T_n(\wt{\phi})+P_nT(\wt{\phi})V_nW_n.
$$
Because $V_n^*=V_{-n}\to0$ strongly, the second term and its adjoint
tend strongly to zero. This settles the last assertion.
\qed

We can now combine the previous result with
Lemma \ref{l1.1} and obtain information about the strong
convergence of the inverses. Note that
$T(\wt{\phi})$ is the transpose of $T(\phi)$.  
\begin{corollary}\label{c1.2}
Let $\phi\in L^\iy(\T)$. If $T_n(\phi)$ is stable, then $T(\phi)$
is invertible and
$$\ba{rcl@{\qquad\qquad}rcl}
T\iv(\phi) &=& \slim\; T_n\iv(\phi), &
T\iv(\wt{\phi}) &=& \slim\; W_nT_n\iv(\phi)W_n.\\
\ea$$
If $M_n(\phi)$ is stable, then $M(\phi)$ and $T(\phi)$
are invertible and
$$\ba{rcl@{\qquad\qquad}rcl}
M\iv(\phi) &=& \slim\; M_n\iv(\phi), &
T\iv(\wt{\phi}) &=& \slim\; W_nM_n\iv(\phi)W_n.\\
\ea$$
\end{corollary}

The next step is the description of invertibility and stability in the
case of functions $\phi$ contained in the set $PC$ of all piecewise
continuous functions on the unit circle.  For this we need more
notation.  In what follows, let $A$ stand for the set $C$ of
continuous functions, the set $\B$, or the set $\Ci$ of infinitely
differentiable functions.  We denote by $\PC{A;K}$ the set of all
functions $\phi$ that can be written in the form
\bq\label{f1.PC}
\phi(e^{i\th}) &=&
b(e^{i\th})t_{\bt_+}(e^{i\th})t_{\bt_-}(e^{i(\th-\pi)})
\prod_{r=1}^R t_{\bt_r^+}(e^{i(\th-\th_r)}) 
t_{\bt_r^-}(e^{i(\th+\th_r)}),
\eq
where $b\in A$ is a nonvanishing function with winding number zero,
$\th_1,\dots,\th_R\in(0,\pi)$ are distinct points, the set
$K\subseteq\{1,-1,e^{i\th_1},\dots,e^{i\th_R},e^{-i\th_1},\dots,
e^{-i\th_R}\}$ is the set of jump discontinuities of $\phi$, and
$\bt_\pm,\bt^\pm_1,\dots,\bt^\pm_R$ are certain complex parameters.
If some these parameters are zero, then no jumps occur at the
corresponding points.

The representation (\ref{f1.PC}) is essentially the same as
(\ref{F.T}), however it displays the special role of jumps at $1$ and
$-1$ and a connection between jumps at a point on the unit circle and
its complex conjugate.  As will be seen below, this distinction is not
necessary for the pure Toeplitz case, however, it is for the Toeplitz
+ Hankel case.  Note that $\PC{C;K}$ is the set of all invertible
functions in $PC$ with finitely many jumps at $K$.

Let $\PCa{A;K}$, $\PCb{A;K}$ and $\PCc{A;K}$ be the sets of functions
of the above form such that in addition the following conditions (I),
(II) and (III), respectively, are satisfied:
\begin{itemize}
\item[(I)]
	$|\Re\bt_\pm|<1/2$ and
	$|\Re\bt_r^\pm|<1/2$ for each $r$;
\item[(II)]
	$-3/4<\Re\bt_{+}<1/4$ and $-1/4<\Re\bt_{-}<3/4$ and
	$|\Re(\bt_r^{+}+\bt_r^{-})|<1/2$ for each $r$;
\item[(III)]
	$-1/2<\Re\bt_{+}<1/4$ and $-1/4<\Re\bt_{-}<1/2$ and
\item[]
	$|\Re\bt_r^+|<1/2$ and $|\Re\bt_r^-|<1/2$ and
	$|\Re(\bt_r^{+}+\bt_r^{-})|<1/2$ for each $r$.
\end{itemize}
We want to emphasize that $\PCc{A;K}\subset\PCa{A;K}\cap\PCb{A;K}$ and
that this inclusion is proper.  This seems strange at first glance,
but it finds its solution in the fact that the representation
(\ref{f1.PC}) and in particular the $\bt$'s need not be unique.
Examples of functions showing this are given in \cite[Sect.~4]{BaEh}.
Finally, note that $\phi\in\PCa{A;K}$ if and only if
$\wt{\phi}\in\PCa{A;\wt{K}}$ where
\bq\label{f.wtK}
\wt{K} &=& \Big\{\;\;1/t\;:\;\;t\in K\;\;\Big\}.
\eq

Invertibility and stability criteria for Toeplitz and Toeplitz +
Hankel operators with piecewise continuous generating functions can
now be stated as follows.  For proofs we refer to \cite{BS} and
\cite{BaEh}, respectively.
\begin{proposition}\label{p1.3}
Let $\phi\in PC$ be a function with jumps in a finite set
$K\subset\T$. Then
\begin{itemize}
\item[(a)] $T(\phi)$ is invertible if and only if $\phi\in\PCa{C;K}$;
\item[(b)] $T_n(\phi)$ is stable if and only if $\phi\in\PCa{C;K}$;
\item[(c)] $M(\phi)$ is invertible if and only if $\phi\in\PCb{C;K}$;
\item[(b)] $M_n(\phi)$ is stable if and only if $\phi\in\PCc{C;K}$.
\end{itemize}
\end{proposition}

We also need to introduce ``approximate'' functions $\phi_\mu$,
$0\le\mu<1$, associated to a piecewise continuous function
$\phi\in\PC{A;K}$ of the form (\ref{f1.PC}):
\bq\label{f1.PCapp}
\phi_\mu(e^{i\th}) &=&
b(e^{i\th})t_{\bt_+,\mu}(e^{i\th})t_{\bt_-,\mu}(e^{i(\th-\pi)})
\prod_{r=1}^R t_{\bt_r^+,\mu}(e^{i(\th-\th_r)}) 
t_{\bt_r^-,\mu}(e^{i(\th+\th_r)})
\eq
Here $t_{\bt,\mu}$ is the smooth nonvanishing function 
with winding number zero defined by
\bq
t_{\bt,\mu}(e^{i\th}) &=& \Big(1-\mu e^{i\th}\Big)^\bt
\Big(1-\mu e^{-i\th}\Big)^{-\bt},\qquad\quad 0\le\th<2\pi.
\eq
Note that also $\phi_\mu\in A$ is a nonvanishing function
with winding number zero.

For a (generalized) sequence $A_\mu$ of operators acting on $\ll$
depending on a parameter $\mu\in[0,1)$, one can define the concepts of
strong convergence (as $\mu\to1$) and stability in the same way as for
(discrete) sequences $A_n$.  The analogues of Lemma \ref{l1.1} and
Lemma \ref{l1.1b} also remain true.  We will write $A=\slim A_\mu$ if both
$A_\mu\to A$ and $A_\mu^*\to A^*$ strongly.
\begin{proposition}\label{p1.4}
Let $K$ be a finite subset of $\T$, let $\phi\in\PC{C;K}$, and
let $\phi_\mu$, $0\le\mu<1$, be the associated approximate functions.
Then
$$
\ba{rcl@{\qquad}rcl@{\qquad}rcl}
H(\phi) &=& \slim H(\phi_\mu), &
T(\phi) &=& \slim T(\phi_\mu), &
M(\phi) &=& \slim M(\phi_\mu).\ea
$$
Moreover, if $\phi\in\PC{\B;K}$ and $f\in\Ci$ vanishes on an open
neighborhood of $K$, then $H(f\phi_\mu)\to H(f\phi)$ and
$H(f/\phi_\mu)\to H(f/\phi)$ in the trace class norm.
\end{proposition}
\proof
First of all note that
$t_\bt(e^{i\th})=(1-e^{i\th})^\bt(1-e^{-i\th})^{-\bt}$, and thus
$$
t_\bt(e^{i\th}) \eql \exp\Big(2i\bt\arg(1-e^{i\th})\Big),\qquad\quad
t_{\bt,\mu}(e^{i\th}) \eql \exp\Big(2i\bt\arg(1-\mu e^{i\th})\Big),
$$
where the argument is taken in $(-\pi/2,\pi/2)$.  It is easy to see
$t_{\bt,\mu}\to t_\bt$ locally uniformly on $\T\setminus\{1\}$ as
$\mu\to1$.  The same holds for the derivatives of arbitrary order.
We obtain that $\phi_\mu\to\phi$ in measure.  Because the sequence
$\phi_\mu$ is uniformly bounded in the norm of $L^\iy(\T)$, the
Laurent operators generated by $\phi_\mu$ (which are unitarily
equivalent to multiplication operators on $L^2(\T)$) converge strongly
to the Laurent operator generated by $\phi$.  This settles the first
part of the proposition.

Now write $\phi_\mu=b\psi_\mu$ and $\phi=b\psi$.  From the above
statements, it follows that $f\psi_\mu\to f\psi$ in the sense of
$\Ci$.  Multiplying with $b\in\B$ we obtain that $f\phi_\mu\to f\phi$
in the norm of $\B$ and hence the desired convergence of the Hankel
operators.  The last assertion can be shown analogously.
\qed

Next we state the necessary and sufficient conditions for the
stability of $T(\phi_\mu)$ and $M(\phi_\mu)$.  Proofs are given in
\cite{ES2} and \cite{BaEh}.  These results can be combined with
Proposition \ref{p1.4} and the analogue of Lemma \ref{l1.1} in order
to obtain a result about the strong convergence of the inverses and
the adjoints of the inverses.  We leave the details to the reader.

\begin{proposition}\label{p1.sconv}
Let $K$ be a finite subset of $\T$, and let $\phi_\mu$, $0\le \mu<1$,
be the approximate functions of the form (\ref{f1.PCapp}) associated
to a function $\phi\in\PC{C;K}$ of the form (\ref{f1.PC}).  Then
\begin{itemize}
\item[(a)]
$T(\phi_\mu)$ is stable if and only if the parameters satisfy
condition (I).
\item[(b)]
$M(\phi_\mu)$ is stable if and only if the parameters satisfy
condition (II).
\end{itemize}
\end{proposition}
Note that it does not suffice to require $\phi\in\PCa{C;K}$ or
$\phi\in\PCb{C;K}$, respectively, because the representation of $\phi$
is not unique and the parameters of $\phi_\mu$ must be chosen
properly.

Finally, we need the following basic results.
\begin{lemma}\label{l1.AB}
Let $A:H_1\to H_2$ and $B:H_2\to H_1$ be linear bounded operators 
acting on
Hilbert spaces $H_1$ and $H_2$. Then the following assertions hold.
\begin{itemize}
\item[(a)]
The operator $I+AB$ is invertible if and only if so is $I+BA$.
\item[(b)]
The operator $I+AB$ is a Fredholm operator if and only if so is 
$I+BA$.
If this is true, then $\ind(I+AB)=\ind(I+BA)$.
\item[(c)]
If $A$ or $B$ is a trace class operator, then $\det(I+AB)=\det(I+BA)$.
\end{itemize}
\end{lemma}
\proof
Part (a) can be proved by using the formula $(I+BA)\iv=I-B(I+AB)\iv
A$.  Assertion (b) can be proved in the same way, by thinking of
the inverses as Fredholm regularizers.  Also, we use the fact that
the kernels (resp.~cokernels) of $I+AB$ and $I+BA$ have the same
dimension.  Finally, $AB$ and $BA$ have the same nonzero eigenvalues
(taking multiplicities into account).
\qed

\begin{lemma}\label{l1.KL}
Let $A_n = P_n+P_nKP_n+W_nLW_n+C_n$ be a sequence of $n\times n$
matrices where $K$ and $L$ are trace class operators, and $C_n$ tends
to zero in the trace class norm.  Then $\lim_{n\to\iy}\det
A_n=\det(I+K)\det(I+L)$.
\end{lemma}
\proof
Because $W_n\to0$ weakly on $\ll$, the sequence $KW_nL$ tends to zero 
in the trace class norm. Hence
\bq
\det A_n &=& \det\Big(I+P_nKP_n+W_nLW_n+C_n\Big)\nn\\
&=& \det\Big((I+P_nKP_n)(I+W_nLW_n)+C_n'\Big)\nn
\eq
with $C_n'\to0$ in the trace class norm. Noting that
$\det(I+W_nLW_n)=\det(I+P_nLP_n)$ completes the proof.
\qed

\section{The limit theorem for piecewise continuous functions}
\label{s.limit}

In Section \ref{s.smooth}, we have proved that for functions
$\phi\in\GB$, the quotient $\det M_n(\phi)/\det T_n(\phi)$ converges
to a nonzero constant.  Surprisingly, the same turns out to be true
for certain piecewise continuous functions satisfying a particular condition
on the location of the jumps.  This condition excludes functions with
jumps at $1$ or $-1$ or at both a point on the unit circle and its
complex conjugate.  The fact that the Hankel operator $H(\phi)$ need
not be compact (and hence $F(\phi)$ need not be defined) makes the
proof of the limit relation for piecewise continuous functions more
complicated than the proof of Proposition \ref{p0.2}.

In order to overcome this obstacle we introduce another operator
determinant which is related to $F(\phi).$  In this
connection, we use the notion of a smooth partition of unity.  By
this we here mean two smooth functions (in $\Ci$) whose sums equals
the constant function with value one on the unit circle.

\begin{theorem}\label{t2.0}
Let $\phi\in L^\iy(\T)$.  Suppose that $T(\phi)$ is invertible and
that there exists a smooth partition of unity, $f+\wt{f}=1$, such
that $H(f\phi)$ and $H(f/\phi)$ are trace class.  Moreover,
introduce the following operators:
\bq
A_{11} &=& T(f)T\iv(\phi)H(\phi),\nn\\[.5ex]
A_{12} &=& T(f)T\iv(\phi)H(\phi)T(\wt{f}),\nn\\[.5ex]
A_{21} &=& T\iv(\phi)H(\phi),\nn\\[.5ex]
A_{22} &=& T\iv(\phi)H(\phi)T(\wt{f}).\nn
\eq
Then $A_{11}$, $A_{12}$ and $A_{22}$ are trace class, and
the operator determinant
\bq
F(\phi;f) &:=& \det \left(
\ba{cc}I+A_{11} & A_{12}\\ -A_{21}A_{11} & I+A_{22}-A_{21}A_{12}\ea
\right)
\eq
is well defined and nonzero.  If in addition $\phi\in\B$, then
$F(\phi;f) = F(\phi)$.
\end{theorem}
\proof
First of all note that if $T(\phi)$ is invertible, then the function
$\phi$ is invertible in $L^\iy(\T)$.  Using the identities
(\ref{f0.Txx}) and (\ref{f0.Hxx}) one can show that $A_{11}$ and
$A_{22}$ are trace class:
\bq
A_{22} &=& T\iv(\phi)H(\phi)T(\wt{f})\nn\\
&=& T\iv(\phi)\Big(H(\phi f)-T(\phi)H(f)\Big)\nn\\
&=& T\iv(\phi)H(\phi f)-H(f),\label{f2.A22}\\[1ex]
A_{11} &=& T(f)T\iv(\phi)H(\phi)\nn\\
&=& 
\Big(T(f/\phi)T(\phi)+H(f/\phi)H(\wt{\phi})\Big)T\iv(\phi)H(\phi)\nn\\
&=& T(f/\phi)H(\phi)+H(f/\phi)H(\wt{\phi})T\iv(\phi)H(\phi)\nn\\
&=& 
H(f)-H(f/\phi)T(\wt{\phi})+H(f/\phi)H(\wt{\phi})T\iv(\phi)H(\phi).
\label{f2.A11}
\eq
Indeed, in each of these terms there appears a Hankel operator which
is trace class.  Hence also $A_{12}$ is trace class.
Notice however that $A_{21}$ need not be trace class.  In any case,
it follows that $F(\phi;f)$ is well defined.  Next, we have
\bq\label{f2.x1}
\left(\ba{cc} I+A_{11} & A_{12} \\ -A_{21}A_{11} & 
I+A_{22}-A_{21}A_{12} \ea\right)
&=&
\left(\ba{cc} I & 0 \\ -A_{21} & I \ea\right)
\left(\ba{cc} I+A_{11} & A_{12} \\ A_{21} & I+A_{22}\ea \right),
\eq
and
\bq\label{f2.x2}
\left(\ba{cc} I+A_{11} & A_{12} \\ A_{21} & I+A_{22} \ea\right)
&=&
\left(\ba{cc} I & 0 \\ 0 & I \ea\right) +
\left(\ba{c} T(f) \\ I \ea\right) T\iv(\phi)H(\phi)
\Big(I\;\;T(\wt{f})\Big).
\eq
Note that $F(\phi;f)\neq0$ if and only if the operator appearing in
(\ref{f2.x1}) is invertible, or, equivalently, if the one in
(\ref{f2.x2}) is invertible.  Using Lemma \ref{l1.AB}(a) and
the fact that
\bq
\Big(I,\;\;T(\wt{f})\Big)\left(\ba{c} T(f) \\ I \ea\right)
&=& I,\nn
\eq
it follows that (\ref{f2.x2}) is invertible if and only if
$I+T\iv(\phi)H(\phi) = T\iv(\phi)M(\phi)$ is invertible.  Hence, we
arrive at the conclusion that $F(\phi;f)\neq0$ if and only if
$M(\phi)$ is invertible.  On the other hand, the operator in
(\ref{f2.x1}) equals identity plus a trace class operator, hence it is
a Fredholm operator with index zero.  It follows that so is
(\ref{f2.x2}).  Using Lemma \ref{l1.AB}(b) we obtain that also
$I+T\iv(\phi)H(\phi)$ and thus $M(\phi)$ is Fredholm with index zero.
Now we need only use the fact that $M(\phi)$ is invertible if and only
if $M(\phi)$ is Fredholm with index zero
(see \cite[Corollary 2.7]{BaEh}).

Finally, suppose that $\phi\in\B$.  Then the operators
$A_{11},A_{12},A_{21},A_{22}$ are trace class, and one can take the
determinant of (\ref{f2.x1}) and (\ref{f2.x2}).  Noting that the
determinant of the first matrix on the right hand side of
(\ref{f2.x1}) is equal to one and employing Lemma \ref{l1.AB}(c), it
follows that $F(\phi;f)=F(\phi)$.
\qed

Although we have not a proof in the general setting, it seems
reasonable that $F(\phi;f)$ does not depend on the particular choice
of $f$ (see also the remark below).

The next result is a quite general version of the limit theorem, which
will be applied afterwards to certain piecewise continuous functions.
\begin{theorem}\label{t2.1}
Let $\phi\in L^\iy(\T)$.  Suppose that the sequence $T_n(\phi)$ is
stable and that there exists a smooth partition of unity,
$f+\wt{f}=1$, such that $H(f\phi)$, $H(f/\phi)$ and
$H(\wt{f}/\wt{\phi})$ are trace class .  Then
$\lim_{n\to\iy}\det M_n(\phi)/\det T_n(\phi)=F(\phi;f)$.
\end{theorem}
\proof
First of all remark that the stability of $T_n(\phi)$ implies the
invertibility of $T(\phi)$.  Hence the assumptions of Theorem
\ref{t2.0} are fulfilled and $F(\phi;f)$ is well defined.  Let
$A_{11},\dots,A_{22}$ be the operators introduced there, and
define the following sequences of matrices:
\bq
A_{11}^{(n)} &=& T_n(f)T_n\iv(\phi)H_n(\phi),\nn\\
A_{12}^{(n)} &=& T_n(f)T_n\iv(\phi)H_n(\phi)T_n(\wt{f}),\nn\\
A_{21}^{(n)} &=& T_n\iv(\phi)H_n(\phi),\nn\\
A_{22}^{(n)} &=& T_n\iv(\phi)H_n(\phi)T_n(\wt{f}).\nn
\eq
We first claim that $A_{11}^{(n)}\to A_{11}$ and $A_{22}^{(n)}\to
A_{22}$ in the trace class norm as $n\to\iy$.  Indeed, using
(\ref{f0.Hxx}), (\ref{f1.VW}) and (\ref{f1.VV}), we can write
\bq
A_{22}^{(n)} &=& T_n\iv(\phi)H_n(\phi)T_n(\wt{f})\nn\\
&=& T_n\iv(\phi)P_nH(\phi)T(\wt{f})P_n-
T_n\iv(\phi)P_nH(\phi)Q_nT(\wt{f})P_n\nn\\
&=& T_n\iv(\phi)P_n\Big(H(\phi f)-T(\phi)H(f)\Big)P_n-
T_n\iv(\phi)P_nV_{-n}H(\phi)H(\wt{f})W_n.\nn
\eq
The last term in the sum tends to zero in the trace class norm because
$H(\phi)H(\wt{f})$ is trace class and $V_{-n}\to0$ strongly.
The first term converges to $A_{22}$ in the trace class
norm since $T_n\iv(\phi)\to T\iv(\phi)$ and $P_n^*=P_n\to I$ strongly
and the expression in the middle is a trace class operator.  Now we
employ identity (\ref{f1.Tn}) to rewrite $T_n(f)$, and we obtain
\bq
A_{11}^{(n)} &=& T_n(f)T_n\iv(\phi)H_n(\phi)\nn\\
&=& T_n(f/\phi)H_n(\phi)+
P_nH(f/\phi)H(\wt{\phi})P_nT_n\iv(\phi)H_n(\phi)\nn\\
&&\mbox{}
+W_nH(\wt{f}/\wt{\phi})H(\phi)W_nT_n\iv(\phi)H_n(\phi).\nn
\eq
Analyzing the first summand yields (see again (\ref{f0.Hxx}),
(\ref{f1.VW}) and (\ref{f1.VV}))
\bq
T_n(f/\phi)H_n(\phi) &=& 
P_nT(f/\phi)H(\phi)P_n-P_nT(f/\phi)Q_nH(\phi)P_n\nn\\
&=& P_n\Big(H(f)-H(f/\phi)T(\wt{\phi})\Big)P_n-W_nH(\wt{f}/\wt{\phi})
H(\phi)V_nP_n.\nn
\eq
Because $H(\wt{f}/\wt{\phi})H(\phi)$ is trace class and
$V_n^*=V_{-n}\to0$ strongly, we obtain that
\bq
T_n(f/\phi)H_n(\phi) &\to& H(f)-H(f/\phi)T(\wt{\phi})\nn
\eq
in the trace class norm as $n\to\iy$. Note that the adjoint of
$T_n\iv(\phi)$ converges strongly to the adjoint of $T\iv(\phi)$.
Hence it is easily seen that
\bq
P_nH(f/\phi)H(\wt{\phi})P_nT_n\iv(\phi)H_n(\phi) &\to&
H(f/\phi)H(\wt{\phi})T\iv(\phi)H(\phi)\nn
\eq
in the trace class norm as $n\to\iy$. Finally, (\ref{f1.VW}) and
(\ref{f1.WW}) imply that
\bq
W_nH(\wt{f}/\wt{\phi})H(\phi)W_nT_n\iv(\phi)H_n(\phi) &=&
W_nH(\wt{f}/\wt{\phi})H(\phi)T_n\iv(\wt{\phi})W_nH_n(\phi)\nn\\&=&
W_nH(\wt{f}/\wt{\phi})H(\phi)T_n\iv(\wt{\phi})T(\wt{\phi})V_nP_n.\nn
\eq
Because $H(\wt{f}/\wt{\phi})H(\phi)$ is trace class, the adjoint of
$T_n\iv(\wt{\phi})$ converges strongly to the adjoint of
$T\iv(\wt{\phi})$ and $V_n^*=V_{-n}\to0$ strongly, it follows that
the latter term converges to zero in the trace class norm.
Thus we have proved that $A_{11}^{(n)}\to A_{11}$ and
$A_{22}^{(n)}\to A_{22}$ in the trace class norm.  Now one can
immediately conclude that also $A_{12}^{(n)}\to A_{12}$
in the trace class norm.
Moreover, it is obvious that $A_{21}^{(n)}\to A_{21}$ strongly.

The desired limit relation can now be proved as follows. Let
\bq
S_n &=& \det M_n(\phi)/\det T_n(\phi) \hsp=\hsp
\det\Big(P_n+T_n\iv(\phi)H_n(\phi)\Big).\nn
\eq
Since $f+\wt{f} = 1$, we have 
\bq
P_n &=& \Big(P_n,\;\;T_n(\wt{f})\Big)\left(\ba{c} T_n(f) \\
P_n\ea\right).\nn
\eq
This in conjunction with Lemma \ref{l1.AB}(c) (for matrices) implies 
that
\bq
S_n &=& \det\left\{P_n+\Big(P_n,\;\;T_n(\wt{f})\Big)\left(
	\ba{c} T_n(f) \\ P_n\ea\right)T_n\iv(\phi)H_n(\phi)\right\}\nn\\
    &=& \det\left\{\left(\ba{cc}P_n & 0 \\ 0 & P_n \ea\right)
	+\left(\ba{c} T_n(f) \\ P_n \ea\right)T_n\iv(\phi)
	H_n(\phi)\Big(P_n,\;\;T_n(\wt{f})\Big)\right\}\nn\\
    &=& \det\left(\ba{cc} P_n + A_{11}^{(n)} & A_{12}^{(n)} \\
        A_{21}^{(n)} & P_n+A_{22}^{(n)} \ea\right)\nn\\
    &=& \det\left(\ba{cc}P_n + A_{11}^{(n)} & A_{12}^{(n)} \\
	-A_{21}^{(n)} A_{11}^{(n)} &
	P_n+A_{22}^{(n)}-A_{21}^{(n)} A_{12}^{(n)}\ea\right)\nn.
\eq
Note that the last identity follows from
\bq
\left(\ba{cc}P_n & 0 \\ -A_{21}^{(n)} & P_n\ea\right)
\left(\ba{cc}P_n + A_{11}^{(n)} & A_{12}^{(n)} \\ A_{21}^{(n)} & 
P_n+A_{22}^{(n)}\ea\right)
&=&\left(\ba{cc}P_n + A_{11}^{(n)} & A_{12}^{(n)} \\
-A_{21}^{(n)} A_{11}^{(n)} &
P_n+A_{22}^{(n)}-A_{21}^{(n)} A_{12}^{(n)}\ea\right),\nn
\eq
where the determinant of the first matrix on the left is equal to one.
Taking the limit $n\to\iy$, we obtain that $S_n \to F(\phi;f)$.
\qed

Note that if the assumptions of Theorem \ref{t2.1} are satisfied
(which are slightly stronger than those of Theorem \ref{t2.0}), then
the value of $F(\phi;f)$ does not depend on the particular choice of
$f$. For it is the limit of a sequence independent of $f$.

Now we specialize to the case of piecewise continuous functions with
jumps on a finite set $K$.  We show how the constant $F(\phi,f)$ can
be evaluated from the constant $F(\phi)$.  The proof will reveal that
the partition of unity condition is in some sense responsible for the
afore--mentioned condition on the location of the jumps of $\phi$.
Note that this condition can be expressed as $K\cap\wt{K}=\emptyset$
where $\wt{K}$ is defined as in (\ref{f.wtK}).

\begin{lemma}\label{l2.5}
Let $\phi\in\PCa{\B;K}$ with $K\cap\wt{K}=\emptyset$.
Moreover, let $\phi_\mu$, $0\le\mu<1$, be the approximate functions
associated to $\phi$.  Then there exists a function $f$ such that the
assumptions of Theorem \ref{t2.0} and Theorem \ref{t2.1} are
fulfilled, and $F(\phi;f) = \lim_{\mu\to 1} F(\phi_\mu)$.
\end{lemma}
\proof
The invertibility of $T(\phi)$ and stability of $T_n(\phi)$ follows
from Proposition \ref{p1.3}(ab).  Because $K\cap\wt{K}=\emptyset$, there
exists an $f\in\Ci$ with $f+\wt{f}=1$ such that $f$ vanishes on an
open neighborhood of $K$.  It is easy to see that $f\phi\in\B$ and
$f/\phi\in\B$.  Thus $H(f\phi)$, $H(f/\phi)$ and $H(\wt{f}/\wt{\phi})$
are trace class.  In addition to $A_{11},\dots,A_{22}$ defined in
Theorem \ref{t2.0}, let
\bq
A_{11}^{(\mu)} &=& T(f)T\iv(\phi_\mu)H(\phi_\mu),\nn\\
A_{12}^{(\mu)} &=& T(f)T\iv(\phi_\mu)H(\phi_\mu)T(\wt{f}),\nn\\
A_{21}^{(\mu)} &=& T\iv(\phi_\mu)H(\phi_\mu),\nn\\
A_{22}^{(\mu)} &=& T\iv(\phi_\mu)H(\phi_\mu)T(\wt{f}).\nn
\eq
Recall that $A_{11}$ and $A_{22}$ can be written in the form
(\ref{f2.A22}) and (\ref{f2.A11}), and analogously so can
$A_{11}^{(\mu)}$ and $A_{22}^{(\mu)}$ with $\phi$ replaced by
$\phi_\mu$.  Applying Proposition \ref{p1.4} and Proposition
\ref{p1.sconv} to these modified expressions, we
conclude that $A_{11}^{(\mu)}\to A_{11}$ and $A_{22}^{(\mu)}\to
A_{22}$ in the trace class norm as $\mu\to 1$.  Moreover, also
$A_{12}^{(\mu)}\to A_{12}$ in the trace class norm and 
$A_{21}^{(\mu)}\to A_{21}$ strongly.  The convergence of these
operators implies that $F(\phi;f)=\lim_{\mu\to1}F(\phi_\mu;f)$.
Finally, observe that $\phi_\mu\in\GB$. Hence
$F(\phi_\mu;f)=F(\phi_\mu)$ by Theorem \ref{t2.0}.
\qed

The evaluation of the constant $F(\phi;f)$ is now
straightforward.  We consider functions $\phi\in\PCa{\B;K}$ with
$K\cap\wt{K}=\emptyset$.  Such a function can be written in the form
(\ref{F.M}).  In terms of this representation, the assumptions on
the parameters of $\phi$ can be expressed as follows:
\begin{itemize}
\item[(i)]
$\bt_+=\bt_-=0$ and $|\Re\bt_r|<1/2$ for each $r$;
\item[(ii)]
$\th_r\in(-\pi,0)\cup(0,\pi)$ for each $r$ and
$\th_r+\th_s\neq0$ for each $r$ and $s$.
\end{itemize}
Hence the final result of this section will nearly treat the situation
which was promised in the introduction.  However, we  still must
exclude jumps at $1$ and $-1$.

\begin{corollary}\label{c2.6}
Let $\phi\in\PCa{\B;K}$ with $K\cap\wt{K}=\emptyset$ be a function
of the form (\ref{F.M}).  Then
\bq
F(\phi;f) &=& F[b]\, \prod_{r=1}^R
	\left(1-e^{-i\th_r}\right)^{\bt_r^2/2+\bt_r/2}
	\left(1+e^{-i\th_r}\right)^{\bt_r^2/2-\bt_r/2}\nn\\
&&\times\prod_{r=1}^R b_+(e^{-i\th_r})^{\bt_r}\,
	\prod_{1\le s<r\le R}
	\Big(1-e^{-i(\th_s+\th_r)}\Big)^{\bt_r\bt_s}.\nn
\eq
\end{corollary}
\proof
Using Theorem \ref{t0.2} and Lemma \ref{l2.5} we can evaluate
$F(\phi_\mu)$ and then pass to the limit $\mu\to1$.
It is easy to see that
\bq
\log\phi_\mu(e^{i\th}) &=& \log b(e^{i\th})+\sum_{r=1}^R
	\bt_r\Big(\log(1-\mu e^{i(\th-\th_r)})-
	\log(1-\mu e^{-i(\th-\th_r)})\Big).\nn
\eq
Employing the Taylor series expansion of $\log(1-z)$ at $z=0$ we obtain
\bq
{[\log \phi_\mu]_n} &=& [\log b]_n -\sum_{r=1}^R \bt_r \mu^n 
e^{-in\th_r}/n,\qquad\quad n>0,\nn\\
{[\phi_{\mu}]_+}(e^{i\th}) &=& 
b_+(e^{i\th})\prod_{r=1}^R \Big(1-\mu 
e^{i(\th-\th_r)}\Big)^{\bt_r}.\nn
\eq
It follows that
\bq
-\frac{1}{2}\sum_{n=1}^\iy n[\log \phi_\mu]_n^2
&=&-\frac{1}{2}\sum_{n=1}^\iy n[\log b]_n^2+\sum_{r=1}^R
\sum_{n=1}^\iy\bt_r [\log b]_n \mu^n e^{-in\th_r}\nn\\ &&
\mbox{}-\frac{1}{2}\sum_{1\le r,s\le R} \sum_{n=1}^\iy
\bt_r\bt_s\mu^{2n}e^{-in(\th_r+\th_s)}/n \nn\\
&=&-\frac{1}{2}\sum_{n=1}^\iy n[\log b]_n^2
+\sum_{r=1}^R \bt_r \log b_+(\mu e^{-i\th_r})\nn\\ &&
\mbox{}+\frac{1}{2}\sum_{1\le r,s\le R}\bt_r\bt_s
\log(1-\mu^2 e^{-i(\th_r+\th_s)}),\nn\\
\left(\frac{\phi_{\mu,+}(1)}{\phi_{\mu,+}(-1)}\right)^{1/2} &=&
\left(\frac{b_+(1)}{b_+(-1)}\right)^{1/2} \prod_{r=1}^R \left(
\frac{1-\mu e^{-i\th_r}}{1+\mu e^{-i\th_r}}\right)^{\bt_r/2}.\nn
\eq
Now it is easy to complete the proof.
\qed

\begin{theorem}\label{t2.7}
Let $\phi\in\PCa{\B;K}$ with $K\cap\wt{K}=\emptyset$ be a function
of the form (\ref{F.M}).  Then
\bq
\det M_n(\phi) &\sim& G[b]^{n}n^{\Omega_M}E_M
\qquad\mbox{ as }n\to\iy.\nn
\eq
\end{theorem}
\proof
The previous results say that
$\det M_n(\phi)/\det T_n(\phi)\to F(\phi;f)$ where 
$F(\phi;f)$ is as given in Corollary \ref{c2.6}. Now we need only
apply the asymptotic formula for Toeplitz determinants
given in (\ref{A.T}).
\qed

\section{The localization theorem for piecewise continuous functions}
\label{s.local}

In this section, we establish a localization theorem for determinants
of Toeplitz + Hankel matrices.  We will show that the quotient $\det
M_n(\phi\psi)/(\det M_n(\phi)\det M_n(\psi))$ converges to some
nonzero constant for certain functions $\phi$ and $\psi$.  Using
this localization, we can reduce the asymptotics of $\det M_n(\phi)$
for certain piecewise continuous functions to the asymptotics for
functions of a particular form with only one or two ``pure'' jumps.

We first define certain constants in terms of operator determinants.
\begin{theorem}\label{t4.EG}
Let $\phi,\psi\in L^\iy(\T)$.
\begin{itemize}
\item[(a)]
Suppose that $T(\phi)$ and $T(\psi)$ are invertible and
$H(\phi)H(\wt{\psi})$ is trace class. Then
\bq
E(\phi,\psi) &=& 
\det\Big(I+T\iv(\phi)H(\phi)H(\wt{\psi})T\iv(\psi)\Big)\nn\\
&=& \det T\iv(\phi)T(\phi\psi)T\iv(\psi)\label{f4.1}
\eq
is well defined and nonzero. Moreover, $E(\wt{\psi},\wt{\phi})$
is well defined and $E(\wt{\psi},\wt{\phi})=E(\phi,\psi)$.
\item[(b)]
Suppose that $M(\phi)$ and $M(\psi)$ are invertible and
$H(\phi)M(\wt{\psi}-\psi)$ is trace class. Then
\bq
G(\phi,\psi) &=& 
\det\Big(I+M\iv(\phi)H(\phi)M(\wt{\psi}-\psi)M\iv(\psi)\Big)\nn\\
&=& \det M\iv(\phi)M(\phi\psi)M\iv(\psi)\label{f4.2}
\eq
is well defined and nonzero.
\item[(c)]
Suppose that $T(\phi)$ and $T(\psi)$ are invertible and 
$H(\phi)H(\wt{\psi})$,
$H(\wt{\phi})H(\psi)$ and $H(\phi)H(\psi)$ are trace class. Then
\bq
H(\phi,\psi) &=& E(\phi,\psi)E(\wt{\phi},\wt{\psi})/E(\phi,\wt{\psi})
\eq
is well defined and nonzero. Moreover, $H(\psi,\phi)$ is well 
defined and
$H(\phi,\psi)=H(\psi,\phi)$.
\end{itemize}
\end{theorem}
\proof
The equality of the expressions in (\ref{f4.1}) and (\ref{f4.2})
follow from (\ref{f0.Txx}) and (\ref{f0.Mxx}).
Note that the operator
\bq\label{f3.h1}
I+T\iv(\phi)H(\phi)H(\wt{\psi})T\iv(\psi) &=&
T\iv(\phi)T(\phi\psi)T\iv(\psi)
\eq
is Fredholm with index zero.  Hence so is $T(\phi\psi)$.  This implies
that $T(\phi\psi)$ is invertible.  It follows that also (\ref{f3.h1})
is invertible and thus $E(\phi,\psi)\neq 0$.  We can argue similarly
in the case of $G(\phi,\psi)\neq 0$.  The point is that the
fact that $M(\phi\psi)$ is Fredholm with index zero implies the
invertibility of $M(\phi\psi)$ (see \cite[Corollary 2.7]{BaEh}).  The
equality $E(\phi,\psi)=E(\wt{\psi},\wt{\phi})$ is obtained by passing
to the transposed operators in (\ref{f4.1}).  
Recall that the transposed operators of $T(\phi)$ and $H(\phi)$
are $T(\wt{\phi})$ and $H(\phi)$, respectively.
Finally, part (c) follows directly from (a).
\qed

Obviously, if $\phi\in\B$, then the above trace class conditions are
fulfilled.  However, the trace class conditions are also fulfilled
under weaker assumptions.
\begin{lemma}\label{l4.trcl}
Let $\phi,\psi\in L^\iy(\T)$.
\begin{itemize}
\item[(a)]
If there exists a smooth partition of unity, $f+g=1$, such that
$f\phi\in\B$ and $g\psi\in\B$, then $H(\phi)H(\wt{\psi})$ is trace 
class.
\item[(b)]
If there exists a smooth partition of unity, $f+g=1$, such that 
$f\phi\in\B$ and $g\wt{\psi}=g\psi$, then 
$H(\phi)M(\wt{\psi}-\psi)$ is trace class.
\end{itemize}
\end{lemma}
\proof
Using equation (\ref{f0.Hxx}) we can write
\bq
H(\phi)H(\wt{\psi})&=&
H(\phi)T(\wt{f})H(\wt{\psi})+H(\phi)T(\wt{g})H(\wt{\psi})\nn\\ 
&=& \Big(H(\phi f)-T(\phi)H(f)\Big)H(\wt{\psi})+H(\phi)
\Big(H(\wt{g\psi})-H(\wt{g})T(\psi)\Big)\label{f4.HH}
\eq
Analogously, from equation (\ref{f0.TM}), we can conclude
\bq
H(\phi)M(\wt{\psi}-\psi) &=&
H(\phi)T(\wt{f})M(\wt{\psi}-\psi)+H(\phi)T(\wt{g})M(\wt{\psi}-\psi)\nn\\ 
&=& \Big(H(\phi f)-T(\phi)H(f)\Big)M(\wt{\psi}-\psi)\nn\\ &&
+H(\phi)\Big(M(\wt{g\psi}-\wt{g}\psi)-H(\wt{g})M(\psi-\wt{\psi})\Big).
\label{f4.HM}
\eq
Hence the operators under consideration are trace class.
\qed

The next result is a first (general) version of a localization theorem
for Toeplitz + Hankel determinants.  All further work is based on it.
\begin{theorem}\label{t4.LocTH}
Let $\phi,\psi\in L^\iy(\T)$ such that the sequences $M_n(\phi)$
and $M_n(\psi)$ are stable.  Suppose in addition that
$H(\phi)M(\wt{\psi}-\psi)$, $H(\phi)H(\psi)$ and
$H(\wt{\phi})H(\psi)$ are trace class.  Then
\bq
\lim_{n\to\iy} \frac{\det M_n(\phi\psi)}{\det M_n(\phi) \det 
M_n(\psi)}
&=& G(\phi,\psi)E(\wt{\phi},\wt{\psi}).\nn
\eq
\end{theorem}
\proof
We obtain from equation (\ref{f0.Mxx}) that
\bq
M_n(\phi\psi) &=& M_n(\phi)M_n(\psi)+P_nH(\phi)M(\wt{\psi}-\psi)P_n
+P_nM(\phi)Q_nM(\psi)P_n,\nn
\eq
and from (\ref{f1.VW}) and (\ref{f1.VV}) it follows that
\bq
P_nM(\phi)Q_nM(\psi)P_n &=&
\Big(P_nT(\phi)V_n+P_nH(\phi)V_n\Big)
\Big(V_{-n}T(\psi)P_n+V_{-n}H(\psi)P_n\Big)\nn\\ 
&=&
W_nH(\wt{\phi})H(\psi)W_n+W_nH(\wt{\phi})H(\psi)V_nP_n\nn\\ &&+
P_nV_{-n}H(\phi)H(\psi)W_n+P_nV_{-n}H(\phi)H(\psi)V_nP_n.\nn
\eq
Because $V_n^*=V_{-n}\to0$ strongly, it is easy to see that the last 
three
terms tend to zero in the trace class norm. Hence
\bq
M_n(\phi\psi) &=& M_n(\phi)M_n(\psi)+P_nH(\phi)M(\wt{\psi}-\psi)P_n
+W_nH(\wt{\phi})H(\psi)W_n+C_n,\nn
\eq
where $C_n\to0$ in the trace class norm.  Multiplying with the
inverses of $M_n(\phi)$ and $M_n(\psi)$ and applying Proposition
\ref{p1.3}(cd) and Corollary \ref{c1.2}, it follows that
\bq
M_n\iv(\phi)M_n(\phi\psi)M_n\iv(\psi) &=&
P_n + 
P_nM_n\iv(\phi)P_nH(\phi)M(\wt{\psi}-\psi)P_nM_n\iv(\psi)P_n\nn\\
&& +W_n\Big(W_nM_n\iv(\phi)W_n\Big)H(\wt{\phi})H(\psi)
   \Big(W_nM_n\iv(\psi)W_n\Big)W_n+C_n'\nn\\
&=&
P_n + P_nM\iv(\phi)H(\phi)M(\wt{\psi}-\psi)M\iv(\psi)P_n\nn\\ &&
+W_nT\iv(\wt{\phi})H(\wt{\phi})H(\psi)T\iv(\wt{\psi})W_n+C_n''\nn
\eq
with $C_n'\to0$ and $C_n''\to0$ in the trace class norm.
Lemma \ref{l1.KL} completes the proof.
\qed

Note that the trace class assumptions required in the theorem can be
replaced by the conditions given in Lemma \ref{l4.trcl}.

Now we proceed with establishing basic properties of the operator
determinants $E(\ast,\ast)$ and $G(\ast,\ast)$ in the case of smooth
functions.  These properties allow their computation.  Note that the
constant $E(\ast,\ast)$ is already known for a long time \cite{B79}.
\begin{theorem}\label{t4.EG2}
Let $b,c\in\GB$.  Then
\bq
E(b,c) &=& F(b)F(\wt{c})/F(b\wt{c}),\label{f4.E2}\\ 
G(b,c) &=& E(b,c)/E(b,\wt{c}).\label{f4.G2}
\eq
In particular, we have
\bq
E(b,c) &=& \exp\left(\sum_{n=1}^\iy n[\log b]_n [\log c]_{-n}\right).
\label{f4.E1}
\eq
\end{theorem}
\proof
We start with proving (\ref{f4.E2}).  Suppose first that $b=\wt{b}$
and $c=\wt{c}$.  Because of (\ref{f0.Mx}) we have $M(bc)=M(b)M(c)$.
Using the definition (\ref{f0.F}) of $F(\ast)$, we can write
\bq
F(b)F(c) &=& \det T\iv(b)M(b)\det T\iv(c)M(c)
\hsp=\hsp \det T\iv(b)M(b)\det M(c)T\iv(c)\nn\\
&=&\det T\iv(b)M(b)M(c)T\iv(c)\hsp=\hsp
\det T\iv(b)M(bc)T\iv(c)\nn\\
&=&\det T\iv(c)T\iv(b)T(bc)\det T\iv(bc)M(bc)
\hsp=\hsp E(b,c)F(bc).\nn
\eq
Hence $F(b)F(c)=E(b,c)F(bc)$.  Now consider arbitrary $b$ and $c$.
Let $b=b_-G[b]b_+$ and $c=c_-G[c]c_+$ be the Wiener--Hopf
factorization.  In the proof of Theorem \ref{t0.2}, we have shown that
$F(b)=F(\wt{b}_+b_+)$.  Analogously, $F(\wt{c})=F(c_-\wt{c}_-)$ and
$F(b\wt{c})=F(\wt{b}_+b_+c_-\wt{c}_-)$.  Doing a similar computation
as in (\ref{f0.TiH}), we obtain
\bq
T\iv(b)H(b)H(\wt{c})T\iv(c) &=&
T\iv(b_+)H(b_+)H(\wt{c}_-)T\iv(c_-).\nn
\eq
Hence $E(b,c) = E(b_+,c_-) = E(\wt{b}_+b_+,c_-\wt{c}_-)$.  Because
$\wt{b}_+b_+$ and $c_-\wt{c}_-$ are even, we can apply the above
results.  It follows
$E(b,c)=E(\wt{b}_+b_+,c_-\wt{c}_-) =
F(\wt{b}_+b_+)F(c_-\wt{c}_-)/F(\wt{b}_+b_+c_-\wt{c}_-) =
F(b)F(\wt{c})/F(b\wt{c})$.  Note that (\ref{f4.E2}) implies
(\ref{f4.E1}) by using Theorem \ref{t0.2}.

We are now going to prove (\ref{f4.G2}). We write
\bq
G(b,c) &=& \det M\iv(b)M(bc)M\iv(c) \nn\\ 
&=& \det M\iv(b)T(b)\det T\iv(b)M(bc)T\iv(c) \det T(c)M\iv(c)\nn\\ 
&=& \det T\iv(c)T\iv(b)M(bc)/(F(b)F(c)) \nn\\ 
&=& \det T\iv(c)T\iv(b)T(bc)\det T\iv(bc)M(bc)/(F(b)F(c))\nn\\ 
&=& E(b,c)F(bc)/(F(b)F(c)) \eql E(b,c)/E(b,\wt{c}).\nn
\eq
Here we have used the equations (\ref{f0.F}), (\ref{f4.1}),
(\ref{f4.2}) and (\ref{f4.E2}).
\qed

Next we address the question under which conditions $E(\ast,\ast)$ and
$G(\ast,\ast)$ are well defined for certain piecewise continuous
functions and how to evaluate them.
\begin{lemma}\label{l4.limit}
Let $K$ and $L$ be finite subsets of $\T$.
\begin{itemize}
\item[(a)]
Let $\phi\in\PCa{\B;K}$ and $\psi\in\PCa{\B;L}$ with $K\cap
L=\emptyset$, and let $\phi_\mu$ and $\psi_\mu$, $0\le\mu<1$, be the
approximate functions associated to $\phi$ and $\psi$.  Then
$E(\phi,\psi)$ is well defined, and
$E(\phi,\psi)=\ds
\lim_{\mu_1\to1}\lim_{\mu_2\to1}E(\phi_{\mu_1},\psi_{\mu_2})=
\lim_{\mu_2\to1}\lim_{\mu_1\to1}E(\phi_{\mu_1},\psi_{\mu_2}).$
\item[(b)]
Let $\phi\in\PCb{\B;K}$ and $\psi\in\PCb{\B;L}$, assume that there
exists an open neighborhood $U$ of $K$ such that
$\psi|_U\equiv\wt{\psi}|_U$, and let $\phi_\mu$ and $\psi_\mu$,
$0\le\mu<1$, be the approximate functions associated to $\phi$ and
$\psi$.  Then $G(\phi,\psi)$ is well defined, and $G(\phi,\psi)=
\ds\lim_{\mu_1\to1}\lim _{\mu_2\to1}G(\phi_{\mu_1},\psi_{\mu_2}).$
\end{itemize}
\end{lemma}
\proof
(a)\ Because $K$ and $L$ are disjoint, there exists a smooth partition
of unity, $f+g=1$, such that $f$ vanishes identically on an open
neighborhood of $K$ and $g$ vanishes on an open neighborhood of $L$.
With this partition the assumptions of Lemma \ref{l4.trcl}(a) are
fulfilled, and hence $H(\phi)H(\wt{\psi})$ is trace class.  We write
this operator and $H(\phi_{\mu_1})H(\wt{\psi}_{\mu_2})$ in the form
(\ref{f4.HH}), and conclude from Proposition \ref{p1.4} that
$H(\phi_{\mu_1})H(\wt{\psi}_{\mu_2})\to H(\phi)H(\wt{\psi})$ in the
trace class norm as $\mu_1\to1$ and $\mu_2\to1$.  Finally, we apply
Proposition \ref{p1.sconv}(a).  Note that the order of
$\mu_1$ and $\mu_2$ in the limit does not play a role.

(b)\ We choose a smooth partition of unity, $f+g=1$, such that $f$
vanishes identically on an open neighborhood of $K$ and $g$ vanishes
on $\T\setminus U$.  Because $f\phi\in\B$ and $g\wt{\psi}=g\psi$,
Lemma \ref{l4.trcl}(b) can be applied, and thus the constant
$G(\phi,\psi)$ is well defined.  As $H(\phi_{\mu_1})$ is trace class,
we conclude from Proposition \ref{p1.4} and Proposition
\ref{p1.sconv}(b) that for fixed $\mu_1$
$$
M\iv(\phi_{\mu_1})H(\phi_{\mu_1})M(\wt{\psi}_{\mu_2}-\psi_{\mu_2})
M\iv(\psi_{\mu_2}) \hsp\to\hsp M\iv(\phi_{\mu_1})H(\phi_{\mu_1})
M(\wt{\psi}-\psi)M\iv(\psi)
$$
in the trace class norm as $\mu_2\to1$.  Writing
$H(\phi_{\mu_1})M(\wt{\psi}-\psi)$ in the form (\ref{f4.HM}), we
obtain from Proposition \ref{p1.4} and Proposition
\ref{p1.sconv}(b) that
$$
M\iv(\phi_{\mu_1})H(\phi_{\mu_1})M(\wt{\psi}-\psi)M\iv(\psi)
\hsp\to\hsp M\iv(\phi)H(\phi)M(\wt{\psi}-\psi)M\iv(\psi)
$$
in the trace class norm as $\mu_1\to1$.  This completes the proof.
\qed

The previous result reduces the evaluation of the constants for
piecewise continuous functions to the case of smooth functions.
However, we must warn the reader to be careful with the choice of
$\phi_\mu$ and $\psi_\mu$.  (see the remark after Proposition
\ref{p1.sconv}).  It may happen, for instance, that given a function
$\phi\in\PCa{\B;K}\cap\PCb{\B;K}$, the associated approximate
functions which satisfy (I) or (II), respectively, are not the same.
Later on, we will restrict to condition (III), and then this
nonuniqueness does not occur.

Now we establish further relations between these constants.
\begin{lemma}\label{l4.EHG}
Let $K$, $L$ and $M$ be finite subsets of $\T$.
\begin{itemize}
\item[(a)]
Let $\phi\in\PCa{\B;K}$, $\psi_1\in\PCa{\B;L}$, $\psi_2\in\PCa{\B;M}$,
and suppose $K\cap(L\cup M)=L\cap M=\emptyset$.
Then $E(\phi,\psi_1\psi_2)=E(\phi,\psi_1)E(\phi,\psi_2)$.
If in addition, $\wt{K}\cap(L\cup M)=\emptyset$, then
$H(\phi,\psi_1\psi_2)=H(\phi,\psi_1)H(\phi,\psi_2)$.
\item[(b)]
Let $\phi\in\PCc{\B;K}$ and $\psi\in\PCc{\B;L}$ with $K\cap
L=K\cap\wt{L}=\emptyset$, and suppose that there exists an open
neighborhood $U$ of $K$ such that $\psi|_U\equiv\wt{\psi}|_U$.
Then $G(\phi,\psi)=E(\phi,\psi)/E(\phi,\wt{\psi})$.
\end{itemize}
\end{lemma}
\proof (a)\
It follows from (\ref{f4.E1}) that the first relation holds for smooth
$\phi$, $\psi_1$ and $\psi_2$.  The general case can is obtained from
Lemma \ref{l4.limit}(a) by approximation.  Note that
$\psi_1\psi_2\in\PCa{\B;L\cup M}$ and
$(\psi_1\psi_2)_\mu=\psi_{1,\mu}\psi_{2,\mu}$ because of
$L\cap M=\emptyset$.  Finally, the relation for $H(\ast,\ast)$ follows
from the definition by using the relation for $E(\ast,\ast)$.

(b)\ 
The equality for smooth functions is stated in (\ref{f4.G2}).  In the
general case, we approximate by smooth functions as indicated in Lemma
\ref{l4.limit}(ab).  Note that the approximate functions $\phi_\mu$
and $\psi_\mu$, respectively, are the same for $G(\ast,\ast)$ and
$E(\ast,\ast)$.
\qed

Now we establish a first version of the localization theorem for the
piecewise continuous functions.
\begin{corollary}\label{c4.LocPC}
Let $\phi\in\PCc{\B;K}$ and $\psi\in\PCc{\B;L}$ with $K$ and $L$ being
finite subsets of $\T$ such that $K\cap L= K\cap\wt{L}=\emptyset$.
Suppose also that there exists an open neighborhood $U$ of $K$ such
that $\psi|_U=\wt{\psi}|_U$.  Then
\bq
\lim_{n\to\iy} \frac{\det M_n(\phi\psi)}{\det M_n(\phi)\det 
M_n(\psi)} &=&
H(\phi,\psi).\nn
\eq
\end{corollary}
\proof
Similar as in the proof of Lemma \ref{l4.limit} (see also Lemma
\ref{l4.trcl}) one can show that the conditions on $K$, $L$ and $U$
imply that the operators $H(\phi)H(\wt{\psi})$, $H(\phi)H(\psi)$ and
$H(\phi)M(\wt{\psi}-\psi)$ are trace class.  Moreover, because of
Proposition \ref{p1.3}, the sequences $M_n(\phi)$ and $M_n(\psi)$ are
stable.  Hence all the assumptions required in Theorem \ref{t4.LocTH}
are fulfilled.  Note that
$H(\phi,\psi)=G(\phi,\psi)E(\wt{\psi},\wt{\phi})$ by Lemma
\ref{l4.EHG}(b).
\qed

The previous result requires a ``strange'' assumption, namely,
$\psi|_U=\wt{\psi}|_U$.  Note that in case $K=\emptyset$
(i.e.~$\phi\in\GB$), one can choose $U=\emptyset$, and hence this
assumption is redundant.  In fact, it turns out that this requirement
can be dropped also in general, as shown in the following
lemma.  For technical reasons we sharpen the smoothness condition.
\begin{lemma}\label{l4.fac}
Let $K$ and $L$ be finite subsets of $\T$ with $K\cap L=K\cap
\wt{L}=\emptyset$, and suppose that $\psi\in\PC{\Ci;L}$.  Then we can
factor $\psi=\psi_1\psi_2$ such that $\psi_1\in C^\iy(\T)$ is
nonvanishing and has winding number zero, and that there exists an
open neighborhood $U$ of $K$ such that
$\psi_2|_U\equiv\wt{\psi}_2|_U$.
\end{lemma}
\proof
Let $U$ and $V$ be open and disjoint neighborhoods of $K$ and $L$,
respectively, which are both non--empty, consist of a finite union of
open subarcs, and satisfy $\wt{U}=U$ and $\wt{V}=V$.  We put
$\psi_1(t)=\psi(t)$ for $t\in U$, and continue $\psi_1$ on
$\T\setminus U$ such that $\psi_1\in\Ci(\T)$ is nonvanishing and has
winding number zero.  This is possible because $\psi$ is infinitely
differentiable and nonzero on $\T\setminus L\supset\T\setminus
V\supset\mbox{clos\,}U$.  The winding number condition can be
fulfilled by choosing the values of $\psi$ appropriately on some
subarc of $V$.  Finally put $\psi_2=\psi/\psi_1$.  Note that
$\psi_2|_U\equiv\wt{\psi}_2|_U\equiv1$.
\qed

\begin{corollary}\label{c4.LocPC2}
Let $\phi\in\PCc{\Ci;K}$ and $\psi\in\PCc{\Ci;L}$ with $K$ and $L$
being finite subsets of $\T$ satisfying
$K\cap L=K\cap\wt{L}=\emptyset$.  Then
\bq
\lim_{n\to\iy} \frac{\det M_n(\phi\psi)}{\det M_n(\phi)\det 
M_n(\psi)} &=& H(\phi,\psi).\nn
\eq
\end{corollary}
\proof
We factor $\psi=\psi_1\psi_2$ as stated in the Lemma \ref{l4.fac}, 
and write
\bq
\frac{\det M_n(\phi\psi)}{\det M_n(\phi)\det M_n(\psi)} &=&
\frac{\det M_n(\phi\psi_1\psi_2)}{\det M_n(\psi_1)
\det M_n(\phi\psi_2)}\times
\frac{\det M_n(\phi\psi_2)}{\det M_n(\phi)\det M_n(\psi_2)}\nn\\ 
&&\times\mbox{}
\frac{\det M_n(\psi_1)\det M_n(\psi_2)}{\det M_n(\psi_1\psi_2)}.\nn
\eq
Because of the conditions on $K$ and $L$ and on the function $\psi_1$,
Proposition \ref{p1.3} shows that the stability of $M_n(\phi)$ and
$M_n(\psi)$ implies the stability of all other sequences $M_n(\ast)$
which occur above.  Employing Corollary \ref{c4.LocPC} one can take
the limit of the above expression, which equals
\bq
\frac{H(\psi_1,\phi\psi_2)H(\phi,\psi_2)}{H(\psi_1,\psi_2)}
&=& H(\psi_1,\phi)H(\phi,\psi_2) \eql H(\phi,\psi).\nn
\eq
Note that Lemma \ref{l4.EHG}(a) and Theorem \ref{t4.EG}(c) has been 
used here.
\qed

\begin{corollary}\label{c4.const}
Let $b\in\GB$,
$\phi(e^{i\th})=t_{\bt_r}(e^{i(\th-\th_r)})$ and
$\psi(e^{i\th})=t_{\bt_s}(e^{i(\th-\th_s)})$. Assume that
$|\Re\bt_r|<1/2$, $|\Re\bt_s|<1/2$, $\th_r,\th_s\in(-\pi,\pi]$,
$\th_r\neq\th_s$ and $\th_r+\th_s\neq0$. Then
\bq
H(b,\psi) &=& b_+(e^{i\th_s})^{\bt_s} b_-(e^{i\th_s})^{-\bt_s}
b_+(e^{-i\th_s})^{\bt_s},\nn\\
H(\phi,\psi) &=& \Big(1-e^{i(\th_s-\th_r)}\Big)^{\bt_r\bt_s}
\Big(1-e^{i(\th_r-\th_s)}\Big)^{\bt_r\bt_s}
\Big(1-e^{-i(\th_r+\th_s)}\Big)^{\bt_r\bt_s}.\nn
\eq
Here $b_\pm$ are the functions defined as in (\ref{2}) and (\ref{3}).
\end{corollary}
\proof
The calculation is similar to the one given in the proof of Corollary
\ref{c2.6}.  We are using (\ref{f4.E1}) and Lemma \ref{l4.limit}(a).
The functions $\phi$ and $\psi$ are approximated by $\phi_\mu$ and
$\psi_\mu$.  We obtain
\bq
E(b,\psi) &=& b_+(e^{i\th_s})^{\bt_s},\nn\\
E(\phi,\psi) &=& \Big(1-e^{i(\th_s-\th_r)}\Big)^{\bt_r\bt_s}.\nn
\eq
Because of $t_\bt(e^{-i\th})=t_{-\bt}(e^{i\th})$, we have
$\wt{\phi}(e^{i\th})=t_{-\bt_r}(e^{i(\th+\th_r)})$ and
$\wt{\psi}(e^{i\th})=t_{-\bt_s}(e^{i(\th+\th_s)})$.
The values of the constant $H(\ast,\ast)$ follow now immediately.
\qed

Now we establish the main result of this section.  We consider
functions $\phi\in\PCc{\B;K}$ of the form (\ref{f1.PC}) and localize
as much as possible.  This allows us to eliminate functions with
pure jumps at $1$ or at $-1$.  However, using this localization
technique one cannot separate singularities at both a point on the unit
circle and its complex conjugate.
\begin{theorem}
Let $\phi$ be a function of the form
\bq
\phi(e^{i\th}) &=&
b(e^{i\th}) \phi^+ (e^{i\th}) \phi^-(e^{i\th})
\prod_{r=1}^R \phi_r(e^{i\th}),\label{f4.LocForm}
\eq
where $b\in\GB$ and
\bq
\phi^+ (e^{i\th}) &=& t_{\bt_+}(e^{i\th}),\nn\\
\phi^-(e^{i\th}) &=& t_{\bt_-}(e^{i(\th-\pi)}),\nn\\
\phi_r(e^{i\th}) &=& t_{\bt^+_r}(e^{i(\th-\th_r)})
t_{\bt^-_r}(e^{i(\th+\th_r)}),\qquad\quad 1\le r\le R.\nn
\eq
Suppose that $\th_1,\dots,\th_R\in(0,\pi)$ are distinct numbers and
that
\begin{itemize}
\item[(a)]$-1/2<\Re\bt_+<1/4$ and $-1/4<\Re\bt_-<1/2$;
\item[(b)]
$|\Re\bt_r^+|<1/2$ and $|\Re\bt_r^-|<1/2$ and
$|\Re(\bt_r^++\bt_r^-)|<1/2$ for each $1\le r\le R$.
\end{itemize}
Then
\bq
\lim_{n\to\iy} \frac{\det M_n(\phi)}{\det M_n(b)\ds
\det M_n(\phi^+ ) \det M_n(\phi^-) \prod_{r=1}^R 
\det M_n(\phi_r)} &=& H,\nn
\eq
where 
\bq
H &=& b_+(1)^{2\bt_+}b_-(1)^{-\bt_+}
b_-(-1)^{2\bt_-}b_-(-1)^{-\bt_-}
2^{3\bt_+\bt_-}\nn\\ &&\times
\prod_{r=1}^R
b_+(e^{i\th_r})^{\bt^+_r+\bt^-_r}
b_-(e^{i\th_r})^{-\bt^+_r}
b_+(e^{-i\th_r})^{\bt^+_r+\bt^-_r}
b_-(e^{-i\th_r})^{-\bt^-_r}
\nn\\ &&\times\prod_{r=1}^R
\Big(1-e^{i\th_r}\Big)^{\bt_+(\bt^+_r+2\bt^-_r)}
\Big(1-e^{-i\th_r}\Big)^{\bt_+(2\bt^+_r+\bt^-_r)}
\nn\\&&\times\prod_{r=1}^R\Big(1+e^{i\th_r}\Big)^{\bt_-(\bt^+_r+2\bt^-_r)}
\Big(1+e^{-i\th_r}\Big)^{\bt_-(2\bt^+_r+\bt^-_r)}
\nn\\ &&\times\prod_{1\le r<s\le R}
\Big(1-e^{i(\th_r+\th_s)}\Big)^{\bt^-_r\bt^-_s+\bt^+_r\bt^-_s+\bt^-_r\bt^+_s}
\Big(1-e^{-i(\th_r+\th_s)}\Big)^{\bt^+_r\bt^+_s+\bt^+_r\bt^-_s+\bt^-_r\bt^+_s}
\nn\\ &&\times\prod_{1\le r<s\le R}
\Big(1-e^{i(\th_r-\th_s)}\Big)^{\bt^+_r\bt^+_s+\bt^-_r\bt^-_s+\bt^-_r\bt^+_s}
\Big(1-e^{-i(\th_r-\th_s)}\Big)^{\bt^+_r\bt^+_s+\bt^-_r\bt^-_s+\bt^+_r\bt^-_s}.\nn
\eq
\end{theorem}
\proof
Note that $\phi$ belongs to $\PCc{\B;K}$ with $K\subseteq
\{1,-1,e^{i\th_1},\dots,e^{i\th_R},e^{-i\th_1},\dots,e^{-i\th_R}\}$
being the set of jump discontinuities of $\phi$, and consequently, so
does each product which involves only some of the factors appearing in
(\ref{f4.LocForm}).  We first apply Corollary \ref{c4.LocPC} and
eliminate the factor $b$.  This yields a constant term
$$
H(b,\phi^+ )H(b,\phi^-)\prod_{r=1}^R H(b,\phi_r)
$$
in the asymptotics (see also Lemma \ref{l4.EHG}(a)).  The remaining
function (i.e.~$\phi$ with $b$ being dropped) is contained even in
$\PCc{\Ci;K}$, and we may 
  
\section{Functions with one jump and the main theorem}
\label{s.1jump}

As noted in previous sections, using the localization technique 
and the limit theorem, we have reduced the computation of the 
asymptotics of the determinants for all piecewise continuous
functions (satisfying appropriate conditions on the size
of the jumps) to those with pure jumps at $1$ or $-1$,
and to those with two jumps at a point on the unit circle and 
its complex conjugate. 

In this section we compute the asymptotics of the corresponding 
determinants for the functions 
$t_{\beta}(e^{i\theta})$ and $t_{\beta}(e^{i(\theta-\pi)})$ and thus
with these examples the promised asymptotic formula given in the 
introduction is proved.

It is interesting to note that we are able to describe the asymptotic
behavior of the determinants for the above pure function
with arbitrary complex parameters $\beta$. Notice that we may exclude
the cases $\beta\in\Z$ as they lead to trivial results.

In computing these examples, and also in the next section where other 
interesting examples are computed, several Cauchy type determinants 
arise. The next lemma shows how to evaluate several of the products 
that occur in the Cauchy determinants in terms of the Barnes 
G-function and how to then evaluate the asymptotics that 
arise from the Barnes G-function. It will be used several times 
in this section and the next.

\begin{lemma}\
(a)\ 
For each nonnegative integer $n$ and each $z\notin\Z$ we have
\bq\label{f3.Bar.form}
\frac{G(1+z-n)}{G(1+z)} &=&
(-1)^{n(n-1)/2}\left(\frac{\sin\pi z}{\pi}\right)^n
\frac{G(1-z+n)}{G(1-z)}.
\eq
(b)\
If $x_1+\dots+x_R=y_1+\dots+y_R$ and
$\omega:=x_1^2+\dots+x_R^2-y_1^2-\dots-y_R^2$, then
\bq\label{f3.Bar.asym}
\prod_{r=1}^R
\frac{G(1+x_r+n)}{G(1+y_r+n)} &\sim&
n^{\omega/2},\qquad\mbox{as } n\to\iy.
\eq
(c)\
For each nonnegative integers $n_1$, $n_2$ and $n$,
the following identities hold:
\bq
\prod_{0\le j<k\le n-1} (k-j) &=& G(1+n),
\label{f3.pr1}\\
\prod_{0\le j<k\le n-1} (k + j +z) &=&\frac{G(2n-1+z)G(\frac{1}{2} 
+\frac{z}{2})G(1+\frac{z}{2})\pi^{\frac{n-1}{2}}}{G(n+z)G(n-\frac{1}{2}
+\frac{z}{2})G(n+\frac{z}{2})2^{(n-1)(n-2 +z)}},
\label{f3.pr4}\\
\prod_{\ts{0\le k_1\le n_1-1 \atop 0\le k_2\le n_2-1}}
(z+k_1+k_2) &=& \frac{G(z+n_1+n_2)G(z)}{G(z+n_1)G(z+n_2)},
\label{f3.pr2}\\
\prod_{\ts{0\le k_1\le n_1-1 \atop 0\le k_2\le n_2-1}}
(z+k_1-k_2) &=& \frac{G(1+z+n_1)G(1-z+n_2)}{G(1+z+n_1-n_2)G(1-z)}
(-1)^{n_2(n_2-1)/2}\left(\frac{\sin\pi z}{\pi}\right)^{n_2}.
\qquad\label{f3.pr3}
\eq
Here we assume $z\notin\{0,-1,-2,\dots\}$ in (\ref{f3.pr4}) and 
(\ref{f3.pr2}), and $z\notin\Z$ in (\ref{f3.pr3}).
\end{lemma}
\proof
(a)\ 
Using the recurrence relation $G(1+z)=\Ga(z)G(z)$ and the  
well known formula $\Ga(z)\Ga(1-z)=\pi/\sin\pi z$, we can write
\bq
\frac{G(1+z-n)}{G(1+z)} &=&
\prod_{k=1}^n \frac{1}{\Ga(1+z-k)} \hsp=\hsp
\prod_{k=1}^n \frac{\sin\pi(k-z)}{\pi}\Ga(k-z) \nn\\ &=&
\prod_{k=1}^n (-1)^{k-1}\frac{\sin\pi z}{\pi}\Ga(k-z) \nn\\ &=&
(-1)^{n(n-1)/2}\left(\frac{\sin\pi z}{\pi}\right)^n \nn
\frac{G(1-z+n)}{G(1-z)}.
\eq

(b)\
In \cite[Proof of Corollary 3.2]{ES} it is shown that
\bq
G(1+z+n) &\sim& a_nb_n^zn^{z^2/2},\qquad\mbox{as }n\to\iy,\nn
\eq
where $a_n$ and $b_n$ are sequences of positive numbers
not depending on $z$.

(c)\
Noting that $G(1)=1$, formula (\ref{f3.pr1})
can be proved as follows:
\bq 
\prod_{j=0}^{n-2}\prod_{k=j+1}^{n-1} (k-j) &=&
\prod_{j=0}^{n-2}\Ga(n-j) 
\eql G(1+n).\nn
\eq
Formula (\ref{f3.pr4}) can be proved in the same sort of way, but 
requires a little more work. First write
\bq
\prod_{0\le j<k\le n-1} (k + j +z) &=& 
\prod_{j=0}^{n-2}\prod_{k=j+1}^{n-1}(j+k+z) \eql
\prod_{j=0}^{n-2}\frac{\Ga(j+n+z)}{\Ga(2j+1+z)}.\nn
\eq
Now we can write this last product as two products, and then 
apply the duplication formula for the Gamma function. We have
\bq
\prod_{j=0}^{n-2}\Ga(j+n+z)\prod_{j=0}^{n-2}\frac{1}{\Ga(2j+1+z)}&=&
\frac{G(2n-1+z)}{G(n+z)}\prod_{j=0}^{n-2}\frac{\pi^{1/2}}{\Ga(j+\frac{1}{2}
+\frac{z}{2})\Ga(j+1+\frac{z}{2})2^{2j+z}},\nn
\eq
and using the basic properties of the Barnes function this is equal to
$$ \frac{G(2n-1+z)G(\frac{1}{2} 
+\frac{z}{2})G(1+\frac{z}{2})\pi^{\frac{n-1}{2}}}{G(n+z)G(n-\frac{1}{2}
+\frac{z}{2})G(n+\frac{z}{2})2^{(n-1)(n-2 +z)}}.$$
Formula (\ref{f3.pr2}) can be shown by writing
\bq
\prod_{k_1=0}^{n_1-1}\prod_{k_2=0}^{n_2-1}(z+k_1+k_2) &=&
\prod_{k_1=0}^{n_1-1}\frac{\Ga(z+n_2+k_1)}{\Ga(z+k_1)} \hsp=\hsp
\frac{G(z+n_1+n_2)G(z)}{G(z+n_1)G(z+n_2)}.\nn
\eq
In order to obtain (\ref{f3.pr3}), we make an index substitution, and
then use (\ref{f3.pr2}):
\bq
\prod_{0\le k_1\le n_1-1 \atop 0\le k_2\le n_2-1} (z+k_1-k_2) &=&
\prod_{0\le k_1\le n_1-1 \atop 0\le k_2\le n_2-1} 
(z+k_1+k_2-n_2+1)\nn\\
&=& \frac{G(1+z+n_1)G(1+z-n_2)}{G(1+z+n_1-n_2)G(1+z)}.\nn
\eq
Finally, we apply (\ref{f3.Bar.form}).
\qed

We now prove the asymptotic formula for the special function 
$t_{\beta}(e^{i\theta}).$ We first note that if $A$ is a matrix
of Cauchy type, that is, if
$\{a_j\}_{j=0}^{n-1}$ and $\{b_k\}_{k=0}^{n-1}$ are sequences of
complex numbers such that the following matrix is well defined
\bq
A &=& \Big[ (a_j+b_k)\iv \Big]_{j,k=0}^{n-1},\nn
\eq
then $\det A =p/q$ where
\bq
p &=& \prod_{0\le j<k\le n-1}(a_k-a_j)(b_k-b_j), \nn\\
q &=& \prod_{0\le j,k\le n-1}(a_j+b_k).\nn
\eq
We will use this identity in the next theorem and also for many of the 
examples that follow.

\begin{theorem}\label{t4a.ex1} 
Let $\phi(e^{i\theta})=t_{\beta}(e^{i\theta})$ and
assume $\beta\notin\Z$.  Then  
\bq
\det M_{n}(\phi) &\sim&   
n^{-3\beta^{2}/2-\beta/2} (2\pi)^{\beta/2} 2^{3\beta^{2}/2}
G(1/2)^{-1} G(1/2-\beta) G(1-\beta) G(1+\beta).\nn
\eq
Moreover, $\det M_n(\phi)=0$ if and only if
$\beta\in\{1/2,3/2,\dots,n-1/2\}$.
\end{theorem}
\proof
The Fourier coefficients of $\phi=t_\beta$ are given by
\bq\label{tFc}
{[t_{\beta}]_n} &=&
\frac{\sin\pi\beta}{\pi(\beta -n)}.
\eq 
Thus the matrices 
$T_{n}(t_{\beta})+H_{n}(t_{\beta})$ have $j,k$ entry 
\bq
\label{f.matr}
\frac{\sin\pi\beta}{\pi}
\left(\frac{1}{\beta -j+k}+\frac{1}{\beta -1-j-k}\right)
&=&
\frac{\sin\pi\beta}{\pi}\cdot
\frac{2\beta-2j-1}{\beta^{2}-\beta-2\beta j+j+j^{2}-k-k^{2}},
\qquad
\eq
and except for terms that can be factored out of rows, 
the corresponding determinant is a Cauchy determinant of the
above form with
$a_j=\beta^{2}-\beta-2\beta j+j+j^{2}$ and
$b_k=-k-k^{2}$.

Our remarks above concerning Cauchy determinants show that
\bq\label{f.term1}
\det M_{n}(t_{\beta}) &=&
\frac{p}{q}\;\left(\frac{\sin\pi\beta}{\pi}\right)^{n}\;
\prod_{j=0}^{n-1}(2\beta -2j -1),
\eq
where
\bq
p &=& \prod_{0\le j<k\le n-1}
(k^{2}+k-2\beta k -j^{2}-j+2\beta j) (-k^{2}-k +j^{2}+j)\nn
\eq
and 
\bq
q &=& \prod_{0\le j,k\le n-1}
(-\beta+j-k)(1-\beta+j+k).\nn
\eq
We can evaluate the $p$ term by first writing 
\bq
p &=& \prod_{0\le j<k\le n-1} (k-j)^{2}(k+j+1)(k+j+1 -2\beta)(-1).\nn
\eq
Then we apply (\ref{f3.pr1}) and (\ref{f3.pr4}) to find
that $p$ is equal to
$$
(-1)^{n(n-1)/2}
\frac{\pi^{n-1} G(1+n)^{2} G(2n) G(1) G(\frac{3}{2})
G(2n-2\beta) G(1-\beta) G(\frac{3}{2}-\beta)}
{2^{2(n-1)(n-1-\beta)}  G(n+1) G(n) G(n+\frac{1}{2})
G(n+1-2\beta) G(n-\beta) G(n+\frac{1}{2}-\beta)}.
$$
To evaluate the $q$ term use (\ref{f3.pr2}) and
(\ref{f3.pr3}) to see that
\bq
q &=& \frac{G(1-\beta+2n) G(1-\beta) G(1-\beta+n)
G(1+\beta+n)}{G(1-\beta+n)^{2} G(1-\beta) G(1+\beta)}
(-1)^{n(n+1)/2}
\left(\frac{\sin\pi\beta}{\pi}\right)^{n}.\nn
\eq
We write the product in (\ref{f.term1}) as
\bq\label{f.term2}
\prod_{j=0}^{n-1}(2\beta -2j-1)  &=&
(-1)^{n}2^{n}\,
\frac{\Ga(n+\frac{1}{2}-\beta)}{\Ga(\frac{1}{2}-\beta)},
\eq
and then simplify and collect terms to obtain
\bq
\det M_{n}(t_{\beta}) &=& 
  \pi^{n-1}  2^{n-2(n-1)(n-1-\beta)}\,
  \frac{\Ga(n+\frac{1}{2}-\beta)}{\Ga(\frac{1}{2}-\beta)}\nn\\
&&\times\;
 \frac{G(n+1)G(2n)G(\frac{3}{2})G(2n-2\beta)G(1-\beta)G(\frac{3}{2}-\beta)}
 {G(n)G(n+\frac{1}{2})G(n+1-2\beta)G(n-\beta)G(n+\frac{1}{2}-\beta)}\nn\\
&&\times\;
 \frac{G(1-\beta+n) G(1+\beta)}
 {G(1-\beta+2n)G(1+\beta+n)},\label{f4.asymfor}
\eq
In the above expression group together the terms with $2n$ factors, 
that is consider
$$
\frac{G(2n)G(2n-2\beta)}{G(1-\beta +2n)}.
$$
We wish to apply  (\ref{f3.Bar.asym}) so we multiply and divide by
$G(2n-\beta-1)$ to find that the terms involving $2n$ are asymptotic to
$$
(2n)^{\beta^{2}-1}G(2n-\beta-1).
$$

Before we evaluate the rest asymptotically it is convenient to use the 
duplication formula for the Barnes G-function \cite{Bar} which reads
\bq\label{Bar.dup}
\textstyle
G(z)G(z+\frac{1}{2})^{2}G(z+1)
&=&\textstyle
G(\frac{1}{2})^{2}  \pi^{z}  2^{(-2z^{2}+3z-1)}  G(2z).
\eq
We let $z=n-\beta/2-1/2$ in this formula to obtain
\bq
G(2n-\beta-1) &=&
  \pi^{-n+\beta/2+1/2}\,
  2^{2n^{2}-2n\beta-5n+\beta^{2}/2+5\beta/2+3}\nn\\
&&\times\;
  \textstyle
  G(n-\frac{\beta}{2}-\frac{1}{2}) G(n-\frac{\beta}{2})^{2}
  G(n-\frac{\beta}{2}+\frac{1}{2}) G(\frac{1}{2})^{-2}\nn
\eq 
and then use this substitution in our formula for
$\det M_{n}(t_{\beta})$. 

So at this point we have, gathering all terms, that
$\det M_{n}(t_{\beta})$ is asymptotically
$$
\begin{array}{c}
\displaystyle
  n^{\beta^{2}-1} \pi^{\beta/2-1/2}\,
  2^{3\beta^{2}/2+\beta/2}\;
  \frac{G(\frac{3}{2}) G(1-\beta)
    G(\frac{3}{2}-\beta) G(1+\beta)}
  {\Ga(\frac{1}{2}-\beta) G(\frac{1}{2})^{2}}\\[2ex]
\displaystyle
\times\;
  \frac{ \Ga(n+\frac{1}{2}-\beta) G(1+n) G(1-\beta+n)
    G(n-\frac{\beta}{2}-\frac{1}{2}) G(n-\frac{\beta}{2})^{2}
    G(n-\frac{\beta}{2}+\frac{1}{2})}
  {G(n) G(n+\frac{1}{2}) G(n+1-2\beta) G(n-\beta)
    G(n+\frac{1}{2}-\beta) G(1+\beta+n)}.
\end{array}
$$
Write
$\Ga(n+\frac{1}{2}-\beta)=G(n+\frac{3}{2}-\beta)/G(n+\frac{1}{2}-\beta)$
and apply (\ref{f3.Bar.asym}) to the above expression in a
straightforward way to finally arrive at
\bq
\det M_{n}(t_{\beta})
&\sim&
  n^{-3\beta^{2}/2-\beta/2}  (2\pi)^{\beta/2}\,
  2^{3\beta^{2}/2}
\textstyle
  G(\frac{1}{2})^{-1}G(1-\beta)G(\frac{1}{2}-\beta)G(1+\beta).\nn 
\eq
This completes the proof of the asymptotic formula. Finally note that the
determinant vanishes if and only if the $p$ term or the product
(\ref{f.term2}) vanish.
\qed
 
For functions that have jump discontinuities at the point $-1$, the 
analogous results of above theorem are contained in the following.  
\begin{theorem}\label{4a.ex2}
Let $\phi(e^{i\theta})=t_{\beta}(e^{i(\theta-\pi)})$
and assume $\beta\notin\Z.$ Then
\bq
\det M_{n}(\phi) &\sim&
n^{-3\beta^{2}/2+\beta/2}
(2\pi)^{\beta/2}  2^{3\beta^{2}/2}
G(3/2)^{-1} G(3/2-\beta)G(1-\beta)G(1+\beta).\nn
\eq
Moreover, $\det M_n(\phi)=0$ if and only if
$\beta\in\{3/2,5/2,\dots,n-1/2\}$.
\end{theorem}
\proof
The only effect of the move of the discontinuity is that we need to 
now evaluate the determinant of $T_{n}(t_{\beta})-H_{n}(t_{\beta})$
since the Fourier coefficients change by a factor of $(-1)^{n}$.
In the computation in the previous theorem this replaces the factor
$(2\beta-2j-1)$ which appears in the numerator of the matrices
(\ref{f.matr}) with $(-2k-1)$.
A simple of check of the computation shows that this 
only changes the first product in the computation, i.e.\
the term $\Ga(n+1/2-\beta)/\Ga(1/2-\beta)$ appearing in (\ref{f.term1})
and in the formulas afterwards has to be replaced by the factor
$\Ga(n+1/2)/\Ga(1/2)$. We leave the details to the reader.
\qed

With these two theorems we have completed all the parts of pieces that 
go together to prove formula (\ref{4}).
For completeness sake we now state the main result 
with all the necessary restrictions on the $\beta$ parameters. The 
theorem follows directly from the two theorems of this section 
combined with the localization result of the last section and 
finally the limit theorem of Section \ref{s.limit}.

\begin{theorem}[Main theorem]
\label{maintheorem}
Let $\phi$ be a function of the form
\bq
\phi(e^{i\th}) &=&
b(e^{i\th}) t_{\beta_+}(e^{i\th})
t_{\beta_-}(e^{i(\th-\pi)})
\prod_{r=1}^R t_{\beta_r}(e^{i(\th-\th_r)}),\nn
\eq
where $b\in\GB$ and where
$\th_1,\dots,\th_R\in(-\pi,0)\cup(0,\pi)$ are distinct numbers 
satisfying $\th_r+\th_s\neq0$ for each $r$ and $s$.
Assume also that
\begin{itemize}
\item[(a)]$-1/2<\Re\bt_+<1/4$ and $-1/4<\Re\bt_-<1/2$;
\item[(b)]$|\Re\bt_r|<1/2$ for each $1\le r\le R$.
\end{itemize}
Then as $n\rightarrow \infty$,
\bq
 \det M_n(\phi) &\sim& G[b]^{n}n^{\Omega_M}E_M,\nn
\eq
where the constants
$G[b]$, $\Omega_M$ and $E_M$ are defined as in
(\ref{5}), (\ref{6}) and (\ref{7}).
\end{theorem}

\section{Other interesting examples}

In this section we evaluate the asymptotics of the
determinant $\det M_n(\phi)$
for some classes of generating functions with two jumps and 
for a class of functions that have a singularity of a different type. 
The functions we consider here are special cases where certain
assumptions on the location and the size of the
singularities are supposed.
Only one class of functions considered here and then only for
certain values of the parameters is covered by the previous theorems. 
We begin by considering the following two functions
$\ph1$ and $\ph2$:
\bq
\ph1(e^{i\th}) &=&\label{phi1}
t_{\bt-1/2}(e^{i\th}) t_{\bt+1/2}(e^{i(\th-\pi)}),\\[1ex]
\ph2(e^{i\th}) &=&\label{phi2}
t_{\bt}(e^{i\th}) t_{\bt}(e^{i(\th-\pi)}).
\eq
These functions have two jumps at $1$ and $-1$. Note that
the functions $\ph1$ with $\bt\in\Z+1/2$ and $\ph2$ with
$\bt\in\Z$ are up to a constant equal to the functions
$\phi(e^{i\th})=e^{in\th}$, $n\in\Z$.
These trivial cases may be excluded without loss of generality.
We remark also that $\ph1$ admits another representation of a similar 
form:
\bq\label{phi1b}
\ph1(e^{i\th}) &=&
-t_{\bt+1/2}(e^{i\th}) t_{\bt-1/2}(e^{i(\th-\pi)}).
\eq

In the next proposition we evaluate the
Fourier coefficients of $\ph1$ and $\ph2$ explicitly.
\begin{proposition}
Let $\ph1$ and $\ph2$ be as above. If $\bt\notin\Z+1/2$, then
\bq
\ph1_n &=& \left\{\ba{ccl}
\ds-\frac{\cos\pi\bt}{\pi(\bt-n/2)}&\quad&
\mbox{if $n$ is odd}\\[2ex]0&&\mbox{if $n$ is even.}\ea\right.
\eq
If $\bt\notin\Z$, then
\bq
\ph2_n &=& \left\{\ba{ccl}
\ds\frac{\sin\pi\bt}{\pi(\bt-n/2)}&\quad&
\mbox{if $n$ is even}\\[2ex]0&&\mbox{if $n$ is odd.}\ea\right.
\eq
\end{proposition}
\proof
Using the definition (\ref{1a}) of the functions $t_{\bt}$, it can be 
easily verified that
\bq
\ph1(e^{i\th}) &=& e^{i\th}t_{\bt-1/2}(e^{2i\th}),\nn\\[1ex]
\ph2(e^{i\th}) &=& t_{\bt}(e^{2i\th}).\nn
\eq
The Fourier series expansion of $\ph1$ and $\ph2$ can be
obtained from that of $t_{\bt-1/2}$ and $t_{\bt}$, respectively.
Recall that the Fourier coefficients of $t_\bt$
has been given in (\ref{tFc}).
\qed

Because of the special form of the Fourier coefficients of $\ph1$ and
$\ph2$, the corresponding matrices $T_n(\phi)+H_n(\phi)$ also have a
particular structure.  In fact, it turns out that they can be
transformed into matrices of Cauchy form.  Thus, using similar
computations from the previous sections the determinants can be
computed.

\begin{proposition}\label{Cauchy2}
Let $\{a_j\}_{j=0}^{m_1-1}$, $\{\ta_j\}_{j=0}^{m_2-1}$, 
$\{b_k\}_{k=0}^{m_1-1}$
and $\{\tb_k\}_{k=0}^{m_2-1}$ such that the block matrix
\bq\label{matrA}
A &=& \left(\ba{cc} A_{11} & A_{12} \\ A_{21} & A_{22} \ea\right)
\eq
is well defined, where
\bq
A_{11} &=& \lb{22ex}{\Big[(a_j+b_k)\iv\Big]_{j,k}}
0\le j\le m_1-1;\;\;\;\; 0\le k\le m_1-1,\nn\\
A_{12} &=& \lb{22ex}{\Big[(a_j+\tb_k)\iv\Big]_{j,k}}
0\le j\le m_1-1;\;\;\;\; 0\le k\le m_2-1,\nn\\
A_{21} &=& \lb{22ex}{\Big[(\ta_j+b_k)\iv\Big]_{j,k}}
0\le j\le m_2-1;\;\;\;\; 0\le k\le m_1-1,\nn\\
A_{22} &=& \lb{22ex}{\Big[(\ta_j+\tb_k)\iv\Big]_{j,k}}
0\le j\le m_2-1;\;\;\;\; 0\le k\le m_2-1.\nn
\eq
Then $\det A=p/q$, where
\bq
p &=& \prod_{0\le j<k\le m_1-1}(a_k-a_j)(b_k-b_j)
\prod_{0\le j<k\le m_2-1}(\ta_k-\ta_j)(\tb_k-\tb_j)%\nn\\&&\times
\prod_{\ts{0\le j\le m_1-1 \atop 0\le k\le 
m_2-1}}(\ta_k-a_j)(\tb_k-b_j),\nn\\[1ex]
q &=& \prod_{0\le j,k\le m_1-1}(a_j+b_k)
\prod_{0\le j,k\le m_2-1}(\ta_j+\tb_k)%\nn\\&&\times
\prod_{\ts{0\le j\le m_1-1 \atop 0\le k\le m_2-1}}(a_j+\tb_k)
\prod_{\ts{0\le j\le m_2-1 \atop 0\le k\le m_1-1}}(\ta_j+b_k).\nn
\eq
\end{proposition}

The matrix $A$ considered in (\ref{matrA}) is also of Cauchy form,
and therefore the above follows from the standard products that arise 
in the Cauchy determinants.
The reason for writing $A$ in block form is only for convenience
in regard to what follows shortly.  

Now we are able to establish the first main results
of this section. In the 
following theorem note that although the values of the parameters
$\beta - 1/2$ and  $\beta +1/2$ never fit the requirements
of the main theorem (Theorem \ref{maintheorem})
for any value of $\beta$,
the answer agrees with the results from that theorem.
This can be seen by a straightforward computation taking into
account the duplication formula for the Barnes G-function
(\ref{Bar.dup}) with $z=1/2-\beta$.
Notice, however, that if we take the
parameters corresponding to representation (\ref{phi1b})
instead of (\ref{phi1}), then the asymptotic formula
given in Theorem \ref{maintheorem} does not
coincide with the following result.

\begin{theorem}\label{t3.th1}
Let $\ph1$ be as defined in (\ref{phi1}) and assume
$\bt\notin\Z+1/2$. Then
\bq
\det M_n(\ph1) &\sim&
n^{-1/4-3\bt^2}2^{4\bt^2}G(1-2\bt)G(1/2+\bt)G(3/2+\bt).\nn
\eq
Moreover, $\det M_n(\ph1)=0$ if and only if 
$\bt\in\{1,2,3,\dots\}$ and $n\geq 2\bt+1$.
\end{theorem}
\proof
For $n$ fixed, put $m_1=m_2=n/2$ if $n$ is even, and put $m_1=(n+1)/2$
and $m_2=(n-1)/2$ if $n$ is odd. Let $\sigma=m_1-m_2$.
Denote by $\ro_{jk}$, $0\le j,k\le n-1$, the
$(j,k)$--entry of the matrix $T_n(\ph1)+H_n(\ph1)$. We permute the 
rows and columns of this matrix in such a way that we take first
the even and then the odd rows and columns.
This rearrangement results in a matrix $B$ with
the same determinant. This matrix is a $2\times2$ block matrix
with the same structure as (\ref{matrA}) and
with a size determined by $m_1$ and $m_2$:
\bq\label{f3.B1}
B &=& \left(\ba{cc} 
\Big[\ro_{2j,2k}\Big]_{j,k}&
\Big[\ro_{2j,2k+1}\Big]_{j,k}\\[2ex]
\Big[\ro_{2j+1,2k}\Big]_{j,k}&
\Big[\ro_{2j+1,2k+1}\Big]_{j,k}\ea\right).
\eq
Because the even Fourier coefficients of $\ph1$ vanish,
it is not hard to see that
\bq\label{f3.B2}
B &=& \left(\ba{cc} 
\Big[\ph1_{2j+2k+1}\Big]_{j,k}&
\Big[\ph1_{2j-2k-1}\Big]_{j,k}\\[2ex]
\Big[\ph1_{2j-2k+1}\Big]_{j,k}&
\Big[\ph1_{2j+2k+3}\Big]_{j,k}\ea\right).
\eq
It follows that
\bq
B &=& -\frac{\cos\pi\bt}{\pi}\left(\ba{cc} I&0\\0&-I\ea\right)A,
\eq
where $A$ is of the form (\ref{matrA}) with
$a_j=\bt-j-1/2$, $\ta_j=-\bt+j+1/2$, $b_k=-k$ and $\tb_k=k+1$.
We obtain that
\bq\label{f3.detB}
\det B &=&
\left(-\frac{\cos\pi\bt}{\pi}\right)^n(-1)^{m_2}\det A,
\eq
where $\det A =p/q$ with
\bq
p &=&
\prod_{0\le j<k\le m_1-1}(k-j)^2
\prod_{0\le j<k\le m_2-1}(k-j)^2
\prod_{\ts{0\le j\le m_1-1 \atop 0\le k\le m_2-1}}
(1-2\bt+k+j)(1+k+j)\nn\\[1ex]
&=& 
G(1+m_1)^2G(1+m_2)^2
\frac{G(1-2\bt+n)G(1-2\bt)}{G(1-2\bt+m_1)G(1-2\bt+m_2)}
\cdot\frac{G(1+n)}{G(1+m_1)G(1+m_2)}\nn\\[1ex]
&=&
\frac{G(1+m_1)G(1+m_2)G(1-2\bt+n)G(1+n)G(1-2\bt)}
{G(1-2\bt+m_1)G(1-2\bt+m_2)}
\label{f3.pp},\\[2ex]
q &=& \prod_{0\le j,k\le m_1-1}(\bt-1/2-j-k)
\prod_{0\le j,k\le m_2-1}(-\bt+3/2+j+k)
\nn\\[1ex]&&\times
\prod_{\ts{0\le j\le m_1-1 \atop 0\le k\le m_2-1}}(\bt+1/2-j+k)
\prod_{\ts{0\le j\le m_2-1 \atop 0\le k\le m_1-1}}(-\bt+1/2+j-k)
\nn\\[1ex]
&=& 
(-1)^{m_1^2}
\frac{G(1/2-\bt+n+\sigma)G(1/2-\bt)}{G(1/2-\bt+m_1)^2}\cdot
\frac{G(3/2-\bt+n-\sigma)G(3/2-\bt)}{G(3/2-\bt+m_2)^2}
\nn\\[1ex]&&\times\;
(-1)^{m_1m_2}
\frac{G(1/2-\bt+m_1)G(3/2+\bt+m_2)}{G(1/2-\bt+\sigma)G(3/2+\bt)}
(-1)^{m_2(m_2-1)/2}\left(-\frac{\cos\pi\bt}{\pi}\right)^{m_2}
\nn\\[1ex]&&\times\;
\frac{G(3/2-\bt+m_2)G(1/2+\bt+m_1)}{G(3/2-\bt-\sigma)G(1/2+\bt)}
(-1)^{m_1(m_1-1)/2}\left(\frac{\cos\pi\bt}{\pi}\right)^{m_1}
\nn\\[1ex]
&=&
(-1)^{n(n-1)/2}\left(-\frac{\cos\pi\bt}{\pi}\right)^n
\frac{G(1/2-\bt+n)G(3/2-\bt+n)}{G(1/2+\bt)G(3/2+\bt)}
\nn\\[1ex]&&\times\;
\frac{G(1/2+\bt+m_1)G(3/2+\bt+m_2)}{G(1/2-\bt+m_1)G(3/2-\bt+m_2)}.
\label{f3.qq}
\eq
In order to obtain these results,
we have first applied Proposition \ref{Cauchy2}, then pulled out a
factor $-1$ in the first and the third product of the $q$--term, 
and finally computed all these products by means of  
(\ref{f3.pr1}), (\ref{f3.pr2}), and (\ref{f3.pr3}).
Recall that $m_1+m_2=n$ and $m_1-m_2=\sigma$, and observe that
$2m_1=n+\sigma$ and $2m_2=n-\sigma$. Note also that
$\sigma=0$ or $\sigma=1$ depending on whether $n$ is even or odd. 
Now we can combine (\ref{f3.detB}), (\ref{f3.pp}) and (\ref{f3.qq}), 
and it follows that
\bq
\det B &=&
G(1-2\bt)G(1/2+\bt)G(3/2+\bt)\cdot 
\frac{G(1-2\bt+n)G(1+n)}{G(1/2-\bt+n)G(3/2-\bt+n)}
\nn\\[1ex]&&\times\;
\frac{G(1+m_1)G(1/2-\bt+m_1)}{G(1-2\bt+m_1)G(1/2+\bt+m_1)}\cdot
\frac{G(1+m_2)G(3/2-\bt+m_2)}{G(1-2\bt+m_2)G(3/2+\bt+m_2)}.\nn
\eq
Here we have used the fact that $n(n-1)/2+m_2$ is always even.
We apply (\ref{f3.Bar.asym}), and we can conclude that the
first fraction in the last expression behaves asymptotically as
$n^{\bt^2-1/4}$, the second
fraction as $(n/2)^{\bt-2\bt^2}$ and the third one as
$(n/2)^{-\bt-2\bt^2}$. This yields the desired limit behavior of
$\det M_n(\ph1)$. It can be read off from the third product in
the $p$--term that the determinant vanishes if and only if
$2\bt\in\{1,\dots,n-1\}$.
\qed
 
The results from the next theorem just as in the previous one agree
with our main theorem (Theorem \ref{maintheorem})
and this example is partially covered by this theorem.
However, we point out that the allowed values of
$\beta_{+}$ and $\beta_{-}$ while of a special form may not satisfy 
the conditions of the main theorem. As before, the duplication
formula for the Barnes G-function shows the equality of both
asymptotic formulas.
 
\begin{theorem}\label{t3.th2}
Let $\ph2$ be as defined in (\ref{phi2}) and assume
$\bt\notin\Z$. Then
\bq
\det M_n(\ph2) &\sim&
n^{-3\bt^2}2^{4\bt^2}G(1-2\bt)G(1+\bt)^2.\nn
\eq
Moreover, $\det M_n(\ph2)=0$ if and only if 
$\bt\in\{1/2,3/2,5/2,\dots\}$ and $n\geq 2\bt+1$.
\end{theorem}
\proof
The proof is similar to the one of the previous theorem. We introduce 
also the numbers $m_1$ and $m_2$ and rearrange the rows and
columns of the matrix $T_n(\ph2)+H_n(\ph2)$ in the same way.
We obtain a matrix $B$ that can be written in the form (\ref{f3.B1}).
However, because now the odd Fourier coefficients vanish,
formula (\ref{f3.B2}) must be modified as follows:
\bq\label{f3.B3}
B &=& \left(\ba{cc} 
\Big[\ph2_{2j-2k}\Big]_{j,k}&
\Big[\ph2_{2j+2k+2}\Big]_{j,k}\\[2ex]
\Big[\ph2_{2j+2k+2}\Big]_{j,k}&
\Big[\ph2_{2j-2k}\Big]_{j,k}\ea\right).
\eq
It follows that
\bq
B &=& \frac{\sin\pi\bt}{\pi}\left(\ba{cc}I&0\\0&-I\ea\right)A,
\eq
where $A$ is of the form (\ref{matrA}) with $a_j=\bt-j$, 
$\ta_j=-\bt+j+1$,
$b_k=k$, $\tb_k=-k-1$. Hence
\bq
\det B &=&
\left(\frac{\sin\pi\bt}{\pi}\right)^n (-1)^{m_2}
\det A,\label{f3.detB2}
\eq
where $\det A=p/q$ with
\bq
p &=& \prod_{0\le j<k\le m_1-1}-(k-j)^2\prod_{0\le j<k\le 
m_2-1}-(k-j)^2
\prod_{\ts{0\le j\le m_1-1 \atop 0\le k\le m_2-1}}
(1-2\bt+k+j)(-1-k-j)\nn\\[1ex]
&=&
(-1)^{n(n-1)/2}G(1+m_1)^2G(1+m_2)^2\nn\\[1ex]&&\times\;
\frac{G(1-2\bt+n)G(1-2\bt)}{G(1-2\bt+m_1)G(1-2\bt+m_2)}
\cdot\frac{G(1+n)}{G(1+m_1)G(1+m_2)}\nn\\[1ex]
&=& (-1)^{n(n-1)/2}
\frac{G(1+m_1)G(1+m_2)G(1-2\bt+n)G(1+n)G(1-2\bt)}
{G(1-2\bt+m_1)G(1-2\bt+m_2)},
\label{f3.ppp}\\[2ex]
q &=&
\prod_{0\le j,k\le m_1-1}(\bt-j+k)
\prod_{0\le j,k\le m_2-1}(-\bt+j-k)\nn\\[1ex]
&&\times\;
\prod_{\ts{0\le j\le m_1-1 \atop 0\le k\le m_2-1}}(\bt-1-j-k)
\prod_{\ts{0\le j\le m_2-1 \atop 0\le k\le m_1-1}}(-\bt+1+j+k)
\nn\\[1ex]
&=&
\frac{G(1+\bt+m_1)G(1-\bt+m_1)}{G(1+\bt)G(1-\bt)}
(-1)^{m_1(m_1-1)/2} \left(\frac{\sin\pi\bt}{\pi}\right)^{m_1}
\nn\\[1ex]&&\times\;
\frac{G(1+\bt+m_2)G(1-\bt+m_2)}{G(1+\bt)G(1-\bt)}
(-1)^{m_2(m_2-2)/2} \left(-\frac{\sin\pi\bt}{\pi}\right)^{m_2}
\nn\\[1ex]&&\times\;
(-1)^{m_1m_2}
\frac{G(1-\bt+n)^2G(1-\bt)^2}{G(1-\bt+m_1)^2G(1-\bt+m_2)^2}
\nn\\[1ex]
&=& (-1)^{n(n-1)/2}(-1)^{m_2}
\left(\frac{\sin\pi\bt}{\pi}\right)^n
\frac{G(1-\bt+n)^2}{G(1+\bt)^2}\cdot
\frac{G(1+\bt+m_1)G(1+\bt+m_2)}{G(1-\bt+m_1)G(1-\bt+m_2)}.
\label{f3.qqq}
\eq
Here we have again employed Proposition \ref{Cauchy2}, pulled out a
factor $-1$ in each of the products of the $p$--term, and finally
evaluated the products by fromula (\ref{f3.Bar.asym}).
Combining (\ref{f3.detB2}), (\ref{f3.ppp}) and (\ref{f3.qqq}) 
it follows that
\bq
\det B &=&
G(1-2\bt)G(1+\bt)^2\cdot
\frac{G(1-2\bt+n)G(1+n)}{G(1-\bt+n)^2}\nn\\[1ex]&&\times\;
\frac{G(1+m_1)G(1-\bt+m_1)}{G(1-2\bt+m_1)G(1+\bt+m_1)}\cdot
\frac{G(1+m_2)G(1-\bt+m_2)}{G(1-2\bt+m_2)G(1+\bt+m_2)}.\nn
\eq
By (\ref{f3.Bar.asym}), the first fraction in this 
expression behaves as $n^{\bt^2}$
and the second and third fraction as $(n/2)^{-2\bt^2}$.
This yields the limit behavior of $\det M_n(\ph2)$.
Finally, the third product in the $p$--term
shows that the determinant vanishes if and only if
$2\bt\in\{1,\dots,n-1\}$.
\qed

Next we consider two further classes of functions $\ph3$ and $\ph4$:
\bq
\ph3(e^{i\th}) &=& t_{\bt-1/2}(e^{i(\th+\pi/2)})
t_{\bt+1/2}(e^{i(\th-\pi/2)}),\\[1ex]
\ph4(e^{i\th}) &=& t_{\bt}(e^{i(\th+\pi/2)})
t_{\bt}(e^{i(\th-\pi/2)}).
\eq
These functions can be directly obtained from $\ph1$ and $\ph2$
if one rotates them by $-\pi/2$ on the unit circle.
The functions $\ph3$ and $\ph4$ have two jumps at $i$ and $-i$.

The following theorem relates the determinants generated by
$\ph3$ and $\ph4$ to those computed in the previous two
theorems.  Note that the corresponding asymptotic formulas,
which can easily be established, cannot be obtained by
piecing together the asymptotics for two functions with a
single jump at $i$ and $-i$.  This fact clearly indicates
the limitations of the localization idea.

\begin{theorem}
Let $\ph1,\dots,\ph4$ be the functions defined above. Then
\bq
\det M_n(\ph3) &=& 
i^{\sigma}\det M_n(\ph1),\\[1ex]
\det M_n(\ph4) &=& \det M_n(\ph2),
\eq
where $\sigma=0$ if $n$ is even and $\sigma=1$ if $n$ is odd.
\end{theorem}
\proof
First note that the Fourier
coefficients are related by $\ph3_n= i^n\ph1_n$.
We make the same rearrangement of the rows and columns
as in the proof of Theorem \ref{t3.th1} and
arrive at a matrix $B$ for $\ph1$ and a matrix
$\wt{B}$ for $\ph3$ both being of the form (\ref{f3.B2}).
We have
\bq
\wt{B} &=& \left(\ba{cc}\diag(i^{2j})_{j=0}^{m_1-1}&0\\
0&\diag(i^{2j})_{j=0}^{m_2-1}\ea\right)B
\left(\ba{cc}\diag(i^{2k+1})_{k=0}^{m_1-1}&0\\
0&\diag(i^{2k-1})_{k=0}^{m_2-1}\ea\right).\nn
\eq
Now we take the determinant.
Analogously, $\ph4_n=i^n\ph2_n$.
After the same modification, we
obtain a matrix $B$ for $\ph2$ and a
matrix $\wt{B}$ for $\ph4$, which are related by
\bq
\wt{B} &=& \left(\ba{cc}\diag(i^{2j})_{j=0}^{m_1-1}&0\\
0&\diag(i^{2j+2})_{j=0}^{m_2-1}\ea\right)B
\left(\ba{cc}\diag(i^{2k})_{k=0}^{m_1-1}&0\\
0&\diag(i^{2k+2})_{k=0}^{m_2-1}\ea\right).\nn
\eq
Taking the determinant completes the proof.
\qed

Among all the functions we have considered so far in this section,
there is only one non-trivial function which is even,
i.e.\ which satisfies $\phi(t)=\phi(1/t)$, $t\in\T$.
This is the function $\ph3$ with $\bt=0$. Note that
\bq
\phi^{(3,0)}(e^{i\th}) &=&
\left\{\ba{rcl} i&\quad&-\pi/2<\th<\pi/2\\[1ex]
-i&&\pi/2<\th<3\pi/2.\ea\right.
\eq

However, there are two more interesting examples that can be done
which are not piecewise continuous functions, but which are even and 
have singularities of Fisher--Hartwig type.

We begin by introducing the function $u_{\alpha},$ defined by
$u_{\alpha}(e^{i\theta}) = (2 -2\cos \theta)^{\alpha}$.
In what follows we assume that $\Re \alpha > -1/2$.
We first note that if we define the functions
\be
\eta_{\gamma}(t) \eql
(1-t)^{\gamma},\qquad\quad
\xi_{\delta}(t) \eql
(1-1/t)^{\delta},
\qquad\quad t\in\T,\label{eta.xi}
\ee
where the branches of $\xi$ and $\eta$ are chosen so that 
$\eta_{\gamma}(0) = \xi_{\delta}(\iy) =1$ for their
analytic continuations, then
$u_{\alpha}= \xi_{\alpha}\eta_{\alpha}$ and
$t_{\beta}= \xi_{-\beta}\eta_{\beta}.$ 
In the following proposition we list some facts
already proven in \cite[Sect.~3]{ES}
which will eventually be used to show how the determinants 
generated by $u_{\alpha}$ 
can be obtained in terms of our previous computations.

\begin{proposition}
Let $D_{\alpha, n}$ be the $n \times n$ diagonal matrix defined by
$$
D_{\alpha, n} \eql \diag
(\mu_{0}^{(\alpha)},\mu_{1}^{(\alpha)},\dots,\mu_{n-1}^{(\alpha)})
$$ 
where
$$\mu_{j}^{(\alpha)} \eql
\frac{\Gamma (1+\alpha+j)}{j!\,\Gamma(1+\alpha)},
\qquad\quad
\Gamma_{\gamma,\delta} \eql 
\frac{\Gamma(1+\gamma)\Gamma(1+\delta)}{\Gamma(1+\delta + \gamma)}.
$$
Then
\bq
T_n(\xi_{\delta}\eta_{\gamma}) &=&
\Gamma_{\delta,\gamma}^{-1} D_{\delta,n}^{-1}
T_n(\eta_{\gamma}) D_{\delta+\gamma,n}
T_n(\xi_{\delta})  D_{\gamma,n}^{-1}\,,\nn
\\[1ex]
H_{n}(\xi_{\delta}\eta_{\gamma}) &=&
\gamma\Gamma_{\delta,\gamma}^{-1} D_{\delta,n}^{-1}
T_n(\eta_{\gamma}) D_{\delta+\gamma,n}
H_n(\tau_{\delta}) D_{-\gamma,n}\,,\nn
\eq
where the finite Hankel matrix 
$H_n(\tau_{\delta})$ is defined by 
\bq
H_n(\tau_{\delta}) &=&
\left(\frac{-(i+j)!\Gamma(1+\delta)}{\Gamma(2+i+j+\delta)}
\right)_{i,j=0}^{n-1}.\nn
\eq
\end{proposition}

The above identities allow us to reduce the computations
of the asymptotics for the function $u_{\alpha}$
to those for a function that we have already done.

\begin{theorem}
Let $\phi(e^{i\theta}) = u_{\alpha}(e^{i\theta})$ with
$\alpha \notin \Z$ and $\Re \alpha > -1/2.$ Then 
\bq
\det M_{n}(\phi) &\sim&
n^{(\alpha^{2}-\alpha )/2} (2\pi)^{-\alpha/2} 2^{3\alpha^2/2}
\frac{G(3/2+\alpha) G(1+\alpha)^{2}}{G(3/2) G(1+2\alpha)}.\nn
\eq
\end{theorem}
\proof
{From} the above identities we may write 
\bq
T_{n}(t_{\beta}) - H_{n}(t_{\beta}) &=&
\Gamma_{-\beta,\beta}^{-1}
D_{-\beta,n}^{-1}  T_{n}(\eta_{\beta})
(T_{n}(\xi_{-\beta}) D_{\beta,n}^{-1} - \beta 
H_{n}(\tau_{-\beta}) D_{-\beta,n})\nn
\eq
and
\bq
T_{n}(u_{\alpha}) + H_{n}(u_{\alpha}) &=&
\Gamma_{\alpha,\alpha}^{-1}
D_{\alpha,n}^{-1}  T_{n}(\eta_{\alpha}) D_{2\alpha,n}
(T_{n}(\xi_{\alpha}) D_{\alpha,n}^{-1} + \alpha
H_{n}(\tau_{\alpha}) D_{-\alpha,n}).\nn
\eq
Pulling out to the right $D_{\beta,n}^{-1}$ and
$D_{\alpha,n}^{-1}$ and taking the determinant, 
this yields with $\beta=-\alpha$,
\bq
\frac{\det(T_{n}(u_{\alpha}) + H_{n}(u_{\alpha}))}
{\det(T_{n}(t_{-\alpha}) - H_{n}(t_{-\alpha}))}
&=&
\frac{\det(\Ga_{\alpha,\alpha}^{-1}
D_{\alpha,n}^{-1}  T_{n}(\eta_{\alpha}) D_{2\alpha,n})}
{\det(\Ga_{\alpha,-\alpha}^{-1}
D_{\alpha,n}^{-1}  T_{n}(\eta_{-\alpha}))}
\cdot
\frac{\det D_{\alpha,n}^{-1}}{\det D_{-\alpha,n}^{-1}}.\nn
\eq
Each of these terms can be computed.  The determinants of 
$T_{n}(\eta_{\alpha})$ and $T_{n}(\eta_{-\alpha})$ are
equal to one since they are triangular matrices.
The other matrices on the right hand side are
diagonal matrices and thus the last term equals
$$
\prod_{j=0}^{n-1}
\frac{\Gamma(1-\alpha+j)\Gamma(1+ 2\alpha+j)}
{j!\, \Gamma(1+\alpha+j)}.
$$
This product can be simplified using the basic recurrence relation
of the Barnes G-function and an application of (\ref {f3.Bar.asym}).
The end result is that asymptotically
\bq\label{f.5asym}
\frac{\det(T_{n}(u_{\alpha}) + H_{n}(u_{\alpha}))}
{\det(T_{n}(t_{-\alpha}) - H_{n}(t_{-\alpha}))}
&\sim&
n^{2\alpha^{2}}\frac{G(1+\alpha)}{G(1-\alpha)G(1+2\alpha)}.
\eq
Finally note that the asymptotics of
$$
\det(T_n(t_{-\alpha})-H_n(t_{-\alpha}))
\eql
\det M_n(t_{-\alpha}(e^{i(\theta-\pi)})),
$$ 
has been computed in Theorem \ref{4a.ex2}.
Collecting all terms gives the desired formula.
\qed

It is easy to modify this last theorem to find the asymptotic formula 
for one last example. We consider $\phi(e^{i\theta}) = 
u_{\alpha}(e^{i(\theta - \pi)}).$ This is a change in the location of 
the singularity/zero to the point $-1.$ As in the examples of the 
previous section this change in the location of the singularity only
requires a small modification in the proof.

\begin{theorem}
Let $\phi(e^{i\theta}) = u_{\alpha}(e^{i(\theta - \pi)})$
with $\alpha\notin\Z$ and $\Re\alpha > -1/2.$ Then 
\bq
\det M_{n}(\phi) &\sim&
n^{(\alpha^{2}+\alpha)/2} (2\pi)^{-\alpha/2} 2^{3\alpha^2/2}
\frac{G(1/2+\alpha) G(1+\alpha)^2}{G(1/2) G(1+2\alpha)}.\nn
\eq
\end{theorem}
\proof
As before the only effect of the move of the discontinuity is that we 
need to evaluate the determinant of $T_{n}(u_{\alpha}) - 
H_{n}(u_{\alpha}).$ This means the only change in the above 
computation is that the term $M_{n}(t_{-\alpha}(e^{i(\theta - \pi)}))$ is 
replaced by $M_{n}(t_{-\alpha}(e^{i\theta})).$ Thus using the results of 
Theorem \ref{t4a.ex1} the asymptotic formula is established.
\qed

Let us make some final remarks concerning the last two theorems.
The assumptions $\alpha\notin\Z$ have been imposed because of
corresponding assumptions in Theorem \ref{t4a.ex1} and
Theorem \ref{4a.ex2}.  However, in the last two theorems
these assumptions are redundant. To see this, one has to elaborate
a bit more on the proofs given in the previous section.
In fact, one can establish an explicit expression for the
determinants. In this expression (using analyticity)
the term $G(1-\alpha)$ appearing in (\ref{f.5asym})
cancels with the term $G(1+\beta)$
appearing in (\ref{f4.asymfor}).

Note that the condition $\Re\alpha>-1/2$ is exactly the condition
for the integrability of the function $u_\alpha$.
Also this condition can be weakened. One can replace it
with the assumption $2\alpha\notin\{-1,-2,\dots\}$.
Notice that in this case one has to
understand $u_{\alpha}$ as a distribution with well
defined Fourier coefficients.
For more details we refer to \cite{ES,Eh}.

\end{document}